\documentclass[12pt]{article}
\usepackage{amssymb}
\usepackage{amsmath, amscd}
\usepackage{color}
\usepackage{latexsym}
\usepackage[dvips]{graphicx}
\usepackage{graphicx}
\usepackage{hyperref}
\usepackage{changes}
\newtheorem{theorem}{Theorem}[section]
\newtheorem{lemma}[theorem]{Lemma}
\newtheorem{corollary}[theorem]{Corollary}
\newtheorem{definition}[theorem]{Definition}
\newtheorem{algorithm}[theorem]{Algorithm}

\newtheorem{remark}{Remark}[section]

\newcommand{\qed}{\hfill $\Box$ }

\newcommand{\proof}{\noindent{\bf Proof.}\ \ }

\usepackage{tikz}
\usetikzlibrary{arrows,automata,positioning}
\usetikzlibrary{graphs,graphs.standard,math}
\usetikzlibrary{shapes.geometric}
\usepackage{textpos}
\usetikzlibrary {quotes}
\usepackage{verbatim}
\usetikzlibrary{arrows,shapes}
\begin{document}

\title{\Large {\bf Isometric embeddings of resonance graphs as finite distributive lattices}}
\author{Zhongyuan Che\\
\normalsize  Department of Mathematics\\
\normalsize  Penn State University, Beaver Campus\\
\normalsize  Monaca, PA 15061, U.S.A.\\
\small {\tt E-mail: zxc10@psu.edu}\\
}

\date{\today}
\maketitle

\begin{abstract}
Let $G$ be a plane bipartite graph and $\mathcal{M}(G)$ be the set of all perfect matchings of  $G$.
The resonance graph $R(G)$ is a graph whose vertex set is $\mathcal{M}(G)$,
and two perfect matchings are adjacent in $R(G)$ if their symmetric difference is a cycle forming the periphery of a finite face of $G$.
It is known that any connected resonance graph can be isometrically embedded  as a finite distributive lattice into hypercubes.  
The isometric dimension of a connected $R(G)$, denoted by $\mathrm{idim}(R(G))$, 
is the smallest dimension of a hypercube that $R(G)$ can be isometrically embedded into. 
Let $d$ be the number of finite faces of $G$ such that there are no forbidden edges on their peripheries.
We show that any connected $R(G)$ has $\mathrm{idim}(R(G)) \ge d$ 
and provide characterizations on when the equality holds. 
Moreover, if a  connected $R(G)$ has $\mathrm{idim}(R(G)) = d$,
then we design an algorithm to generate a binary coding of length $d$  for all perfect matchings of $G$ which induces an
isometric embedding of $R(G)$  as a finite distributive lattice into a $d$-dimensional hypercube without generating $\mathcal{M}(G)$.
Our results provide answers for the fundamental cases of  both open questions 
raised in [\textit{SIAM J. Discrete Math.} {\bf 22} (2008) 971--984.]

\vskip 0.2in
\noindent {\emph{keywords}}:  finite distributive lattice, forcing face, isometric dimension, peripheral convex expansion,   
plane (weakly) elementary bipartite graph, resonance (di)graph

\end{abstract}

\section{Introduction} 
A \textit{perfect matching} (or, a \textit{1-factor}) of a graph is a set of vertex disjoint edges 
incident to all vertices of the graph.
A \textit{benzenoid graph} (or, a \textit{benzenoid system}) 
is a 2-connected plane bipartite graph such that each finite face is a unit hexagon.
A benzenoid graph is called \textit{catacondensed} if all vertices are on the periphery of the graph.
A perfect matching is called a \textit{Kekul\'e structure} of a benzenoid graph in chemistry. 
An edge of a graph with a perfect matching is called \textit{forbidden} if it is not contained in any perfect matching of the graph.
A bipartite graph is \textit{elementary} if and only if it is connected and without any forbidden edges \cite{LP86}. 
A plane bipartite graph (not necessarily connected) with a perfect matching is \textit{weakly elementary} 
if the removal all forbidden edges  of the graph does not produce any new finite faces,
and the resulted components are called \textit{elementary components}.
 By definitions, we can see that any plane elementary bipartite graph is also weakly elementary.
It is well known \cite{Z06} that benzenoid graphs are plane weakly elementary bipartite graphs.
The infinite face of a plane bipartite graph $G$ is \textit{forcing} if the periphery of $G$  is an even cycle and the 
subgraph obtained by removing all vertices on the periphery of $G$ is either empty or has exactly one perfect matching.
Therefore, catacondensed benzenoid graphs and 2-connected outerplane bipartite graphs are
plane elementary bipartite graphs whose infinite face is forcing.

The concept of resonance graphs was first introduced by chemists 
for benzenoid graphs and closely related to the resonance energy of benzenoids.
It was extended to plane bipartite graphs in terms of \textit{perfect matching graphs}  
and \textit{$Z$-transformation graphs} by mathematicians. 
An early survey paper on resonance graphs was given by Zhang in 2006 \cite{Z06}, and
more structural properties of resonance graphs have been obtained since then. 
For example, see \cite{BCTZ23+, BCTZ24+, C18, C19, C21, Z10, ZLS08} and their reference papers.  
The resonance graph of a plane bipartite graph $G$ is connected if and only if $G$ is weakly elementary  \cite{F03, ZZY04},
and any connected resonance graph is a median graph \cite{ZLS08}, which is a special type of partial cubes.
Hence, any connected resonance graph can be isometrically embedded into hypercubes. 
A partial order $\le_L$ can be defined on the set of all perfect matchings $\mathcal{M}(G)$ of a plane bipartite graph $G$ 
to form a finite distributive lattice $\mathbb{M}(G)=(\mathcal{M}(G), \le_L)$ if and only if $G$ is weakly elementary \cite{Z10},
and the Hasse diagram of  $\mathbb{M}(G)$ is isomorphic to a digraph of $R(G)$ \cite{LZ03}.
Hence, any connected resonance graph can be oriented as the Hasse diagram of a finite distributive lattice.    

It was shown in \cite{ZLS08} that if $G$ is a plane weakly elementary bipartite graph with $n$ finite faces,
then there is a nonnegative integer coding of length $n$ on the finite distributive lattice 
$\mathbb{M}(G)$ preserving its partial order.
Moreover,  such a nonnegative integer coding can be transformed into 
a binary coding of length $d$  which is order-preserving on  the finite distributive lattice $\mathbb{M}(G)$ 
and distance-preserving on the resonance graph $R(G)$, where $d$ is equal to both the height of $\mathbb{M}(G)$ and the diameter of $R(G)$.
For the special case when $G$ is a catacondensed benzenoid graph or a 2-connected outerplane bipartite graph, 
such a binary coding has length $n$,
and can be generated for all perfect matchings of $G$ without generating $\mathcal{M}(G)$.
The motivation of this paper is mainly from the two open problems raised in  \cite{ZLS08}.

\begin{itemize}
\item {Open problem 1.} Determine the smallest integer $m$ such that $\mathbb{M}(G)$ can be embedded into $(\mathbb{Z}^+)^m$, 
where $\mathbb{Z}^+$ is the set of nonnegative integers.

\item {Open problem 2.} For a general benzenoid system $G$, how can one design an algorithm to generate a binary coding for all 1-factors of $G$,
avoiding the generation of $\mathcal{M}(G)$?
\end{itemize}
 
In this paper, we first focus on the case when $G$ 
is a plane elementary bipartite graph with $n$ finite faces for some positive integer $n$.
We show that $\mathrm{idim}(R(G)) \ge n$ and the equality holds if and only if the infinite face of $G$ is forcing.
We design an algorithm to generate a binary coding of length $n$ for all perfect matchings of $G$ whose 
infinite face is forcing without knowing $\mathcal{M}(G)$, and 
such a binary coding induces an isometric embedding of $R(G)$ as a finite distributive lattice into 
a hypercube $Q_n$.
This is a generalization of Algorithm 5.1 and Theorem 5.2 given by Zhang et al. in \cite{ZLS08} 
from resonance graphs of catacondensed benzenoid graphs and 2-connected outer plane bipartite graphs 
to resonance graphs of plane elementary bipartite graphs whose infinite face is forcing.
In general, if $G$ is a plane bipartite graph whose resonance graph $R(G)$ is connected, and $d$ is 
the number of finite faces of $G$ without forbidden edges on their peripheries,
then $\mathrm{idim}(R(G)) \ge d$ and the equality holds if and only if 
$G$ is a plane weakly elementary bipartite graph whose each elementary component
with more than two vertices has the property that the infinite face is forcing,
and a binary coding for all perfect matchings of $G$ can be generated 
to induce an isometric embedding of $R(G)$  as a finite distributive lattice into $Q_d$
without generating $\mathcal{M}(G)$.
Our main results solve the fundamental cases of Open problem 1 and Open problem 2 raised in \cite{ZLS08}.
Improved versions with $\mathrm{idim}(R(G))$ are also obtained for Theorem 4.3 and Corollary 4.4 in \cite{ZLS08}.

\section{Preliminaries}\label{S:Pre}

\subsection{Partial cubes and median graphs} 

Let $G$ be a graph with the vertex set $V(G)$ and the edge set $E(G)$. 
The distance between two vertices $u$ and $v$ of $G$
is the length of a shortest path between them in $G$, and denoted by $d_G(u,v)$.
An induced subgraph $H$ of $G$ is called an \textit{isometric subgraph}   if 
for any two vertices $u$ and $v$ of $H$, $d_H(u,v)=d_G(u,v)$.
\textit{Partial cubes} are isometric subgraphs of hypercubes.
The \textit{isometric dimension} of a partial cube $G$, denoted by $\mathrm{idim}(G)$, is
the least integer $n$ for which $G$ can be isometrically embedded into an $n$-dimensional hypercube \cite{D73}.
Two edges $x_1x_2$ and $y_1y_2$ of a connected graph $G$ are  
in \textit{relation $\Theta$} if $d_G(x_1, y_1) + d_G(x_2, y_2) \neq  d_G(x_1, y_2) + d_G(x_2, y_1)$.
Relation $\Theta$ is reflexive and symmetric but not necessarily transitive. 
For any partial cube $G$, relation $\Theta$ is transitive and forms an equivalence relation on $E(G)$,
and the number of $\Theta$-classes of $G$ is the isometric dimension of $G$ \cite{D73}.

The set of all shortest paths between two vertices $u$ and $v$ in a graph $G$ is denoted by $I_G(u,v)$.
A connected graph $G$ is called a \textit{median graph} if $I_G(x, y) \cap I_G(y, z) \cap I_G(x, z)$ has exactly one vertex
for any three vertices $x, y, z$ of $G$. Any median graph is a partial cube \cite{HIK11}, 
and so can be isometrically embedded into hypercubes.
A connected subgraph $H$ of $G$ is \textit{convex} in $G$ if for any two vertices $x$ and $y$ of $H$,
 $I_G(x,y) \subseteq I_H(x,y)$. A well known structure characterization of a median graph is stated as a sequence of convex expansions.
Let $G_1$ and $G_2$ be isometric subgraphs of $G$ such that $V(G)=V(G_1) \cup V(G_2)$.
Let $G_1 \cap G_2$ be an induced subgraph of $G$ generated by  $V(G_1) \cap V(G_2)$. 
Assume that $G_1 \cap G_2$ is a nonempty subgraph of $G$, 
and  there are no edges between $V(G_1) \setminus V(G_2)$ and $V(G_2)\setminus V(G_1)$.
Take disjoint copies of $G_1$ and $G_2$, 
and connect every vertex of $G_1 \cap G_2$ in $G_1$ 
with the same vertex of $G_1 \cap G_2$ in $G_2$ with an edge. 
Then the resulted graph is called an \textit{expansion} of $G$.
In particular, if $G_1=G$,
then the above expansion is called a \textit{peripheral expansion} of $G$.
An expansion of $G$ is a \textit{convex expansion}  if $G_1 \cap G_2$ is a convex  subgraph of $G$.
A  graph is a median graph if and only if it can be obtained from the one-vertex graph 
by a sequence of convex expansions \cite{HIK11}. 
This implies that each edge of a median graph is contained in an induced 4-cycle of the graph.

\subsection{A finite distributive structure of resonance graphs}

Let $G$ be a graph with a perfect matching. A subgraph $H$ of $G$ is called \textit{nice} if there is a perfect matching
on the subgraph of $G$ obtained by removing all vertices $H$ from $G$.
A cycle (or, a path) of $G$ with a perfect matching $M$ is called \textit{$M$-alternating} 
if its edges are in and  out of $M$ alternately.  Let  $M_1$ and $M_2$ be  two distinct perfect matchings of $G$.
The symmetric difference of $M_1$ and $M_2$, denoted by $M_1 \oplus M_2$, 
is the set of edges contained in either $M_1$ or $M_2$ but not both. 
A cycle (or, path)  of $G$ is  \textit{$(M_1, M_2)$-alternating} if its edges are in $M_1$ and $M_2$ alternately.
It is well known \cite{LP86} that $M_1 \oplus M_2$ is a set of vertex disjoint cycles 
that are $(M_1,M_2)$-alternating. 

All plane bipartite graphs considered in this paper have perfect matchings unless specified otherwise.
Let $G$ be a plane bipartite graph and $\mathcal{M}(G)$ be the set of all perfect matchings of $G$. 
Let  $M \in \mathcal{M}(G)$. 
A face of $G$ is called \textit{$M$-resonant} if the periphery of the face is an $M$-alternating cycle,
and we say that a face is \textit{resonant} briefly if there is no need to specify the perfect matching $M$.
It is known \cite{ZZ00} that every face (including the infinite face) of a plane elementary bipartite graph with more than two vertices is resonant.
Assume that all vertices of $G$ are colored black and white  such that adjacent vertices cannot have the same color.
Then an $M$-alternating cycle $C$ of $G$ is called \textit{$M$-proper} (respectively,  \textit{$M$-improper}) 
if every edge of $C$ contained in $M$ goes from 
white to black vertices (respectively,  from black to white vertices) along the clockwise orientation of $C$.
A plane bipartite graph $G$ with a perfect matching has a unique perfect matching $M_{\hat{0}}$ (respectively,  $M_{\hat{1}}$)
such that $G$ has no proper $M_{\hat{0}}$-alternating cycles (respectively,  no improper $M_{\hat{1}}$-alternating cycles) \cite{ZZ97}.

The {\em resonance graph}  of  a plane bipartite graph $G$, denoted by $R(G)$,  is the graph 
whose vertices are all perfect matchings of $G$, and there is an edge
between two perfect matchings $M_1,M_2$ of $G$ 
if the symmetric difference $M_1 \oplus M_2$ is the facial cycle of a finite face $s$ of $G$ \cite{ZZY04}, 
and we say that the edge $M_1M_2$ of $R(G)$ has the \textit{face-label} $s$.
The  \textit{resonance digraph} (or, \textit{$Z$-transformation digraph}), denoted by $\overrightarrow{R}(G)$,
is the digraph obtained from $R(G)$ by adding a direction for each edge  so that
$\overrightarrow{M_1M_2}$ is an directed edge from $M_1$ to $M_2$ if $M_1 \oplus M_2$ 
is a proper $M_1$-alternating (or, an improper $M_2$-alternating) facial cycle for a finite face of $G$.
Let $\mathcal{M}(G)$ be the set of all perfect matchings of $G$. 
Then a partial order $\le_L$ can be defined on $\mathcal{M}(G)$ such that $M' \le_L M$ 
if there is a directed path from $M$ to $M'$ in $\overrightarrow{R}(G)$.

A \textit{poset}  $\mathbb{S}$ is a set  $S$ with a partial order $\le$.
Let $x$ and $y$ be two distinct elements of $S$.
We say that $y$ \textit{covers} $x$ if $x \le y$ and there is no other element $z \in S \setminus \{x, y\}$
such that $x  \le z \le y$.
The \textit{Hasse diagram} of $\mathbb{S}$ is a digraph with the vertex set $S$ such that 
there is a directed edge from $y$ to $x$ if and only if $y$ covers $x$.  
A \textit{lattice} is a poset that is closed under two operations \textit{join} and \textit{meet},
where the \textit{join} of $x$ and $y$ is the unique least upper bound of $x$ and $y$,
and  the \textit{meet} of $x$ and $y$ is the unique greatest lower bound of $x$ and $y$.
A lattice is \textit{distributive} if two operations join and meet admit distributive laws, 
and  \textit{finite} if it has finite many elements.
Every finite lattice has a maximum element and a minimum element.
The \textit{height} of a finite lattice is the length of a directed path from the maximum element to the minimum element
in its Hasse diagram.

Let $G$ be a plane weakly elementary bipartite graph. 
Two functions $\phi_{M_1}$ and $\psi_{M_1M_2}$ where $M_1, M_2 \in \mathcal{M}(G)$ were introduced in \cite{LZ03}
to show that $\mathbb{M}(G)=(\mathcal{M}(G), \le_L)$ is a finite distributive lattice  whose 
Hasse diagram is isomorphic $\overrightarrow{R}(G)$.
Two functions $\phi_{M_1}$ and $\psi_{M_1M_2}$  were also applied in  \cite{ZLS08} to show that 
$R(G)$ is a median graph,  
and obtain the related results about nonnegative integer codings and binary codings on the finite distributive lattice 
$\mathbb{M}(G)$.

\begin{definition}\label{D:Functions}\cite{LZ03}
Let $G$ be a plane weakly elementary bipartite graph with $n$ finite faces for some positive integer $n$.
Let $\mathcal{M}(G)$ be the set of all perfect matchings of $G$.
Assume that $M_1, M_2 \in \mathcal{M}(G)$, and $f$ is a finite face of $G$.
Let $\phi_{M_1}(f)$  be the number of cycles in $M_1 \oplus M_{\hat{0}}$ 
with $f$ in their interiors, and such cycles are all proper $M_1$-alternating and pairwise vertex disjoint if they exist.
Let $\psi_{M_1M_2} (f)$ be the number of proper $M_1$-alternating cycles in $M_1 \oplus M_2$
with $f$ in their interiors minus the number of improper $M_1$-alternating cycles in $M_1 \oplus M_2$
with $f$ in their interiors. 
\end{definition}
 
The following lemma will be needed to obtain our main results.

\begin{lemma} \label{L:phi-psi} \cite{ZLS08} 
Let $G$ be a plane weakly elementary bipartite graph with $n$ finite faces  for some positive integer $n$.
Let $\mathcal{M}(G)$ be the set of all perfect matchings of $G$.
Then for any $M_1, M_2 \in \mathcal{M}(G)$, $\phi_{M_1}-\phi_{M_2}=\psi_{M_1M_2}$.
Moreover,  if $s_1, s_2, \ldots, s_n$ are finite faces of $G$, then $M_1$ covers $M_2$ if and only if 
$M_1 \oplus M_2$ is the facial cycle of $s_i$ for some $1 \le i \le n$  such that $\phi_{M_1}(s_i)-\phi_{M_2}(s_i)=1$
and $\phi_{M_1}(s_j)-\phi_{M_2}(s_j)=0$  for all other finite faces $s_j$ where $1 \le j \neq i \le n$.
\end{lemma}

\vskip 0.1in

\subsection{A decomposition structure of resonance graphs}

A cycle that is the periphery of a face in a plane bipartite graph $G$ is called a \textit{facial cycle}.
A facial cycle of a finite face $s$ is denoted by $\partial s$. The facial cycle of the infinite face
is the periphery of $G$ and can be denoted by $\partial G$. 

\begin{definition}\label{D:Paths}
Let $P$ be an odd length path in a graph $G$ with a perfect matching. 
Let $\mathcal{M}(G)$ be the set of all perfect matchings of $G$.
Then two subsets of $\mathcal{M}(G)$ can be defined.

Let $\mathcal{M}(G; P^{-})$ be the set of all perfect matchings $M$ of $G$ 
such that $P$ is $M$-alternating and without end edges contained in $M$.

Let $\mathcal{M}(G; P^{+})$ be the set of all perfect matchings $M$ of $G$ 
such that $P$ is $M$-alternating and with end edges contained in $M$.

Furthermore, if  $P$ is the odd length path of a plane bipartite graph $G$ and located on the periphery of a finite face $s$ of $G$,
then two more subsets of $\mathcal{M}(G)$ can be defined.

Let $\mathcal{M}(G; P^{-}, \partial s)$ be the set of all perfect matchings $M$ of $G$ 
contained in $\mathcal{M}(G; P^{-})$ such that $s$ is $M$-resonant, that is, $\partial s$
is $M$-alternating.

Let $\mathcal{M}(G; P^{+}, \partial s)$ be the set of all perfect matchings $M$ of $G$ 
contained in $\mathcal{M}(G; P^{+})$ such that $s$ is $M$-resonant,  that is, $\partial s$
is $M$-alternating.
\end{definition}

Let $G$ be a graph with a perfect matching. 
If $P$ is an odd length path of $G$, then $\mathcal{M}(G; P^{-})$ and $\mathcal{M}(G; P^{+})$ do not necessarily form a partition of $\mathcal{M}(G)$.
A \textit{handle}  of a graph is a path that has exactly two end vertices with degree larger than 2 in the graph \cite{CC13}. 
If $P$ is an odd length handle of $G$, then $\mathcal{M}(G; P^{-})$ and $\mathcal{M}(G; P^{+})$ form  a partition of $\mathcal{M}(G)$.
The concept of a handle was introduced in \cite{C18} as a key tool to partition the set of all perfect matchings of a plane bipartite graph
and obtain a decomposition structure of its resonance graph.

Let $G$ be a plane elementary bipartite graph. A finite face $s$ of $G$  is called a \textit{reducible face} of $G$ if the
common periphery of $s$ and $G$ is a handle $P$ of odd length, and
the removal of the internal vertices and edges of $P$ results in a plane elementary bipartite graph. 
Based on the facts that $R(G)$ is a median graph \cite{ZLS08} and that any median graph has a decomposition structure with respect
to a $\Theta$-class \cite{HIK11},
we showed that the set of all edges whose face-label is a reducible face of $G$ is a $\Theta$-class of $R(G)$,
and provided a  decomposition structure of the resonance graph $R(G)$ with respect to a reducible face of $G$ in \cite{C18}.

Let $X$ be a subset of $\mathcal{M}(G)$, which is the set of all perfect matchings of a plane bipartite graph $G$. 
We use $\langle X \rangle$ to denote the induced subgraph of the resonance graph $R(G)$ generated by $X$.
 
\begin{theorem}\label{T:MedianR(G)}\cite{C18}
Assume that $G$ is a plane elementary bipartite graph and $s$ is a reducible face of $G$. 
Let $P$ be the common periphery of $s$ and $G$.
Let $H$ be the subgraph of $G$ obtained by removing all internal vertices and edges of $P$.
Let $F$ be the set of all edges in the resonance graph $R(G)$ with the face-label $s$. 
Then $F$ is a $\Theta$-class of $R(G)$ and $R(G)-F$ has exactly two components
$\langle \mathcal{M}(G; P^{-}) \rangle$ and $\langle \mathcal{M}(G; P^{+}) \rangle$, 
where $\langle \mathcal{M}(G; P^{-})\rangle $ is isomorphic to the resonance graph $R(H)$,
and the following properties hold true: (i) $F$ is a matching defining an isomorphism between 
$\langle \mathcal{M}(G; P^{-}, \partial s) \rangle$ and $\langle \mathcal{M}(G; P^{+}, \partial s) \rangle$;
(ii) $\langle \mathcal{M}(G; P^{-}, \partial s) \rangle$ is convex in $\langle \mathcal{M}(G; P^{-}) \rangle$,
and $\langle \mathcal{M}(G; P^{+}, \partial s) \rangle$ is convex in $\langle \mathcal{M}(G; P^{+}) \rangle$;
(iii) Both $\langle \mathcal{M}(G; P^{-}) \rangle$ and $\langle \mathcal{M}(G; P^{+}) \rangle$ are median graphs.
Furthermore, $R(G)$ can be obtained from $R(H)$ by a peripheral convex expansion 
if and only if $\mathcal{M}(G; P^{+})=\mathcal{M}(G; P^{+}, \partial s)$.
\end{theorem} 

\begin{corollary}\label{C:MedianR(G)}\cite{C18}
Assume that $G$ is a plane elementary bipartite graph and $s$ is a reducible face of $G$. 
Let $P$ be the common periphery of $s$ and $G$.
Let $H$ be the subgraph of $G$ obtained by removing all internal vertices and edges of $P$.
If $s$ has common edges with exactly one other finite face of $G$,
then the resonance graph $R(G)$ can be obtained from the resonance graph $R(H)$ by a peripheral convex expansion.
\end{corollary}

A decomposition construction of a special type of plane bipartite graphs was introduced in \cite{ZZ00}.
Start from an edge $e$, add a path $P_1$ of odd length which has two end vertices in common with these of $e$
such that $e$ and $P_1$ form an even cycle $G_1$ which is the periphery of a finite face $s_1$,
proceed inductively to build a sequence of plane bipartite
graphs  $G_i$ for $2 \le i \le n$ where $G_n=G$  as follows.
If $G_{i-1} =e + P_1+ \ldots + P_{i-1}$ has already been constructed, 
then $G_i=G_{i-1}+P_i$ can be obtained by adding the $i$th path $P_i$ of odd length
in the exterior of $G_{i-1}$ 
such that $P_i$ has exactly two end vertices in common with $G_{i-1}$,
$P_i$ and a part of the periphery of $G_{i-1}$ surround a finite face $s_i$.
By definition, we can see that each $G_i$ is a plane elementary bipartite graph, and $s_i$ is a reducible face of $G_i$
such that the common periphery of $s_i$ and $G_i$ is the path $P_i$ for all $2 \le i \le n$.
The above construction is called a \textit{reducible face decomposition} of $G$, 
and denoted briefly as  $\mathrm{RFD}(G_1, G_2, \ldots, G_n)$.

It was shown \cite{ZZ00} that a plane bipartite graph with more than two vertices is elementary 
if and only if it has a reducible face decomposition. 
Note that the resonance graph of an even cycle graph is the one-edge graph. 
Therefore, by Theorem \ref{T:MedianR(G)}, if $G$ is a plane elementary bipartite graph with more than two vertices,
then its resonance graph $R(G)$ can be constructed from the one-edge graph 
by a sequence of convex expansions with respect to a reducible face decomposition of $G$.
By Corollary \ref{C:MedianR(G)}, for a special case that $G$ has a $\mathrm{RFD}(G_1, G_2, \ldots, G_n)$ 
with the property that the reducible face $s_i$ of $G_i$ has common edges 
with exactly one other finite face $s_{\alpha(i)}$ in $G_i$ for all $2 \le i \le n$,
then $R(G)$ can be constructed from the one-edge graph 
by a sequence of peripheral convex expansions
with respect to a reducible face decomposition of $G$.

\begin{remark} If $G$ is a plane elementary bipartite graph whose $R(G)$ can be constructed from the one-edge graph 
by a sequence of peripheral convex expansions
with respect to a reducible face decomposition of $G$,
then it is not necessary that $G$ has a $\mathrm{RFD}(G_1, G_2, \ldots, G_n)$ 
with the property that the reducible face $s_i$ of $G_i$ has common edges 
with exactly one other finite face $s_{\alpha(i)}$ in $G_i$ for all $2 \le i \le n$.
For example, see Figure 1 and Figure 2 in \cite{C21}.
\end{remark}

\subsection{When the infinite face is forcing}

It was shown in \cite{CC13} that any 2-connected plane bipartite graph with a forcing face (including the infinite face) is elementary.
We organize related known results with references or proofs if needed for completeness and clarity.

\begin{lemma}\label{L:Basic}
Let $G$ be a plane elementary bipartite graph with $n$ finite faces for some positive integer $n$. 
Let $\mathcal{M}(G)$ be the set of all perfect matchings of $G$, and
 $\phi_M$ be the function defined in Definition \ref{D:Functions} for any perfect matching $M \in \mathcal{M}(G)$.
Then the following statements are equivalent.

(i) $\phi: M \rightarrow \phi_M$ defines a binary coding of length $n$  
on $\mathcal{M}(G)$ which
induces an isometric embedding of $R(G)$  as a finite distributive lattice into a hypercube $Q_n$ 
whose Hasse diagram is isomorphic to $\overrightarrow{R}(G)$.

(ii) Any two vertex disjoint cycles forming a nice subgraph of $G$ have disjoint interiors.

(iii) The infinite face of $G$ is forcing.

(iv)  $R(G)$ can be constructed from the one-edge graph 
by a sequence of peripheral convex expansions with respect to a reducible face decomposition of $G$.

\end{lemma}
\proof By Theorem 4.4 in \cite{LZ03}, $\mathbb{M}(G)=(\mathcal{M}(G), \le_L)$ is a finite distributive lattice
whose Hasse diagram is isomorphic to $\overrightarrow{R}(G)$.
By \cite{LZ03, ZLS08}, $\phi: M \rightarrow \phi_M$ defines a nonnegative integer coding 
$\phi_M=\phi_M(s_1)\phi_M(s_2) \ldots \phi_M(s_n)$ of length $n$
for each $M \in \mathcal{M}(G)$ which is order-preserving on the finite distributive lattice $\mathbb{M}(G)$.
Furthermore, by Theorem 4.3 in \cite{ZLS08},  $\phi: M \rightarrow \phi_M$  is a binary coding of length $n$ 
if and only if any two vertex disjoint cycles of $G$ forming a nice subgraph have disjoint interiors. 
Hence, $(i)$ and $(ii)$ are equivalent.

By Lemma 3.2 in \cite{BCTZ24+},  $(ii)$ and $(iii)$ are equivalent. 

By Theorem 2 in \cite{C21}, $(iii)$ and $(iv)$ are equivalent.
 \qed\\

\section{Main Results}
\subsection{When $G$ is a plane elementary bipartite graph}  

\begin{theorem}\label{T:IsometricDimension}
Let $G$ be a plane elementary bipartite graph with $n$ finite faces for some positive integer $n$. 
Then $\mathrm{idim}(R(G)) \ge n$ and the equality holds if
and only if the infinite face of $G$ is forcing.
Moreover, if $\mathrm{idim}(R(G)) = n$, then $\phi: M \rightarrow \phi_M$ 
defines a binary coding of length $n$ for all perfect matchings of $G$ which
induces an isometric embedding of $R(G)$ into a hypercube $Q_n$ as a finite distributive lattice 
whose Hasse diagram is isomorphic to $\overrightarrow{R}(G)$.
\end{theorem}
\proof  By \cite{D73}, the isometric dimension of a partial cube is  the number of $\Theta$-classes of the partial cube.
By  \cite{HIK11}, any two edges on a shortest path of graph cannot be in the same $\Theta$-class.
Hence, the isometric dimension of a partial cube is at least the diameter of the partial cube.
By Theorem 3.1 in \cite{ZLS08}, $R(G)$ is a median graph,  so $R(G)$ is a partial cube.
Hence, $\mathrm{idim}(R(G))$ is at least the diameter of $R(G)$.

By Corollary 3.5 and Theorem 4.6 in \cite{ZLS08}, 
there is a binary coding of length $d$ on the finite distributive lattice $\mathbb{M}(G)=(\mathcal{M}(G), \le_L)$
which is order-preserving and induces an isometric embedding of $R(G)$ into the hypercube $Q_d$,
 where $d$ is equal to both the height of $\mathbb{M}(G)$ and the diameter of $R(G)$.
It follows that $\mathrm{idim}(R(G))$ is both the height of $\mathbb{M}(G)$ and the diameter of $R(G)$.

By the proof of Theorem 2 in \cite{C21}, the height of $\mathbb{M}(G)$ is at least $n$, and the equality holds
if and only if  $R(G)$ can be obtained from 
the one-edge graph by a sequence of peripheral convex expansions with respect to a reducible face decomposition of $G$
if and only if the infinite face of $G$ is forcing.
Therefore, $\mathrm{idim}(R(G)) \ge n$ and the equality holds if and only if the infinite face of $G$ is forcing.
The conclusion follows by Lemma \ref{L:Basic}. \qed\\

\begin{corollary} \label{C:RFD-Special} 
Let $G$ be a plane elementary bipartite graph with a $\mathrm{RFD}(G_1, G_2, \ldots, G_n)$ 
associated with a sequence of finite faces $s_1, s_2, \ldots, s_n$ for some positive integer $n$
such that $s_i$ is a reducible face of $G_i$ for all $2 \le i \le n$. 
If each $s_i$ has common edges with exactly one other finite face $s_{\alpha(i)}$ in $G_i$ for $2 \le i \le n$,
then $\mathrm{idim}(R(G))=n$ and $\phi: M \rightarrow \phi_M$ defines a binary coding of length $n$ for all perfect matchings of $G$ which
induces an isometric embedding of $R(G)$ into a hypercube $Q_n$ as a finite distributive lattice 
whose Hasse diagram is isomorphic to $\overrightarrow{R}(G)$. 
\end{corollary}
\proof 
Note that $G_1$ is an even cycle and $R(G_1)$ is the one-edge graph. 
Applying Corollary \ref{C:MedianR(G)}  repeatedly for $G_i$ with a reducible face $s_i$ 
which has common edges with exactly one other finite face $s_{\alpha(i)}$ in $G_i$ for $2 \le i \le n$.
Then $R(G)$ can be obtained from 
the one-edge graph by a sequence of peripheral convex expansions with respect to a reducible face decomposition of $G$.
By Lemma \ref{L:Basic}, the infinite face of $G$ is forcing. Hence, the conclusion follows by Theorem \ref{T:IsometricDimension}.
\qed\\

The following algorithm is a generalization of Algorithm 5.1 in \cite{ZLS08} 
for catacondensed benzenoid graphs and 2-connected outerplane bipartite graphs.

\begin{algorithm}\label{A:FDL}
{\bf Input:} A plane elementary bipartite graph $G$ whose infinite face is forcing.  
Assume that $G$ has a $\mathrm{RFD}(G_1, G_2, \ldots, G_n)$   
associated with a sequence of finite faces $s_1, s_2, \ldots, s_n$ for some positive integer $n$
such that $s_i$ is a reducible face of $G_i$ for all $2 \le i \le n$. 

{\bf Output:} A binary coding list $\mathbf{L_n}$  for all perfect matchings of $G$ which
induces an isometric embedding of $R(G)$ into a hypercube $Q_n$ as a finite distributive lattice 
whose Hasse diagram is isomorphic to $\overrightarrow{R}(G)$.

Step 0. $i := 1$, $\mathbf{L_i} := \{0,1\}$ 

Step 1. If $i = n$, stop.

Step 2. 
Let $P_{i+1}$ be the common periphery of $s_{i+1}$ and $G_{i+1}$.
If $P_{i+1}$  has two end vertices colored white and black respectively 
along the clockwise orientation of the periphery of $G_{i+1}$, 
set $\mathbf{L_{i+1}} := \{x0 : x \in \mathbf{L_i}\}  \cup  \{x1 : x \in \mathbf{L_i}$ and $(x)_{j} =1$ 
for all $j$ such that $s_j$ is a finite face of $G_{i}$  having  common edges with $s_{i+1} \}$;
Otherwise, set $\mathbf{L_{i+1}} := \{x1 : x \in \mathbf{L_i}\}  \cup  \{x0 : x \in \mathbf{L_i} $ and $(x)_{j} =0$ 
for all $j$ such that $s_j$ is a finite face of $G_{i}$ having  common edges with $s_{i+1} \}$,
where $(x)_{j}$ is the $j$-th digit of a binary string $x$ of length $i$.

Step 3. $i:=i+1$, go to step 1.
\end{algorithm}

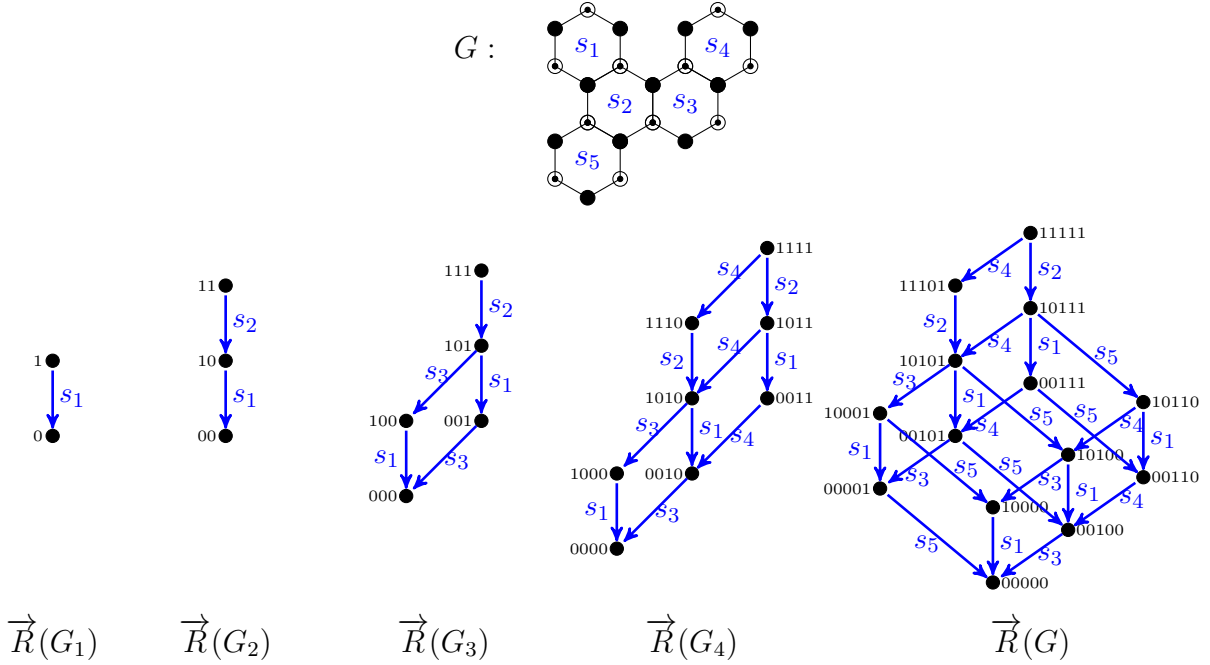
\begin{figure}
\begin{center}
\tikzset{Bullet/.style={draw=black, circle, minimum size=3pt, scale=0.5}}
\begin{tikzpicture}[scale=0.5]
  \begin{scope}[xshift=-3cm]
   \node [black]{$G:$};
  \end{scope}
 
    \draw [color=black, mark=*] plot[samples at={-150,-90,...,150,210},  variable=\x] 
  (\x:1);
    \node [Bullet] at (-150:1){};
  \node  [Bullet, fill] at (-90:1) {};
  \node [Bullet] at (-30:1) {} ;
  \node [Bullet, fill] at (30:1) {} ;
  \node  [Bullet] at (90:1) {};
  \node [Bullet, fill] at (150:1) {};
  \node [blue]{$s_1$};

 \begin{scope}[xshift=0.866 cm,  yshift=-1.5cm]
    \draw [color=black,mark=*] plot[samples at={-150,-90,...,210},variable=\x] 
  (\x:1);
      \node [Bullet] at (-150:1){};
  \node  [Bullet, fill] at (-90:1) {};
  \node [Bullet] at (-30:1) {} ;
  \node [Bullet, fill] at (30:1) {} ;
  \node  [Bullet] at (90:1) {};
  \node [Bullet, fill] at (150:1) {};
  \node[blue]{$s_2$};
  \end{scope}
    \begin{scope}[xshift=2.598 cm,  yshift=-1.5cm]
    \draw [color=black,mark=*] plot[samples at={-150,-90,...,210},variable=\x] 
  (\x:1);
      \node [Bullet] at (-150:1){};
  \node  [Bullet, fill] at (-90:1) {};
  \node [Bullet] at (-30:1) {} ;
  \node [Bullet, fill] at (30:1) {} ;
  \node  [Bullet] at (90:1) {};
  \node [Bullet, fill] at (150:1) {};
   \node [blue]{$s_3$};
  \end{scope}
      \begin{scope}[xshift=3.464 cm]
    \draw [color=black, mark=*] plot[samples at={-150,-90,...,210},variable=\x] 
  (\x:1);
      \node [Bullet] at (-150:1){};
  \node  [Bullet, fill] at (-90:1) {};
  \node [Bullet] at (-30:1) {} ;
  \node [Bullet, fill] at (30:1) {} ;
  \node  [Bullet] at (90:1) {};
  \node [Bullet, fill] at (150:1) {};
   \node [blue]{$s_4$};
  \end{scope}
  \begin{scope}[yshift=-3cm]
    \draw [color=black,mark=*] plot[samples at={-150,-90,...,210},variable=\x] 
  (\x:1);
      \node [Bullet] at (-150:1){};
  \node  [Bullet, fill] at (-90:1) {};
  \node [Bullet] at (-30:1) {} ;
  \node [Bullet, fill] at (30:1) {} ;
  \node  [Bullet] at (90:1) {};
  \node [Bullet, fill] at (150:1) {};
   \node [blue]{$s_5$};
  \end{scope}
\end{tikzpicture}

\tikzstyle{every node}=[draw,shape=circle, inner sep=0pt, minimum size=5pt, label distance=0.01pt, node distance=1cm]
\begin{tikzpicture}[>=stealth', shorten > = 1pt,  el/.style = {inner sep=2pt, align=left},every label/.append style = {font=\tiny}]
\node [fill, label=left:$0$] (1) {};
\node [fill, label=left:$1$] (2) [above of=1] {};
\path[every node/.style={font=\sffamily\small}, inner sep=0.2pt, line width=1.pt]
[<-] (1) edge [color=blue] node[el, right, pos=0.5] {$s_1$}  (2);
\node[draw=none, below of =1, yshift=-1.7cm] {$\overrightarrow{R}(G_1)$};

\begin{scope}[xshift=2.3cm]
\node [fill, label=left:$00$] (1) {};
\node [fill, label=left:$10$] (2) [above of=1] {};
\node [fill, label=left:$11$] (3) [above of=2] {};
\path[every node/.style={font=\sffamily\small}, inner sep=0.5pt, line width=1.pt]
[<-] (1) edge [color=blue] node[el, right, pos=0.5] {$s_1$} (2)
[<-] (2) edge [color=blue] node[el, right, pos=0.5] {$s_2$} (3);
\node[draw=none, below of =1, yshift=-1.7cm]{$\overrightarrow{R}(G_2)$};
\end{scope}

\begin{scope}[xshift=5.7cm, yshift=0.2cm]
\node [fill, label=left:$001$] (1) {};
\node [fill, label=left:$101$] (2) [above of=1] {};
\node [fill, label=left:$111$] (3) [above of=2] {};
\node [fill, label=left:$000$] (4) [below of=1, xshift=-1cm] {};
\node [fill, label=left:$100$] (5) [below of=2, xshift=-1cm] {};
\path[every node/.style={font=\sffamily\small}, inner sep=0.5pt, line width=1.pt]
[<-] (1) edge [color=blue] node[el, right, pos=0.5] {$s_1$} (2)
[<-] (4) edge [color=blue] node[el, left, pos=0.5] {$s_1$} (5)
[<-] (2) edge [color=blue] node[el, right, pos=0.5] {$s_2$} (3)
[->] (1) edge [color=blue ] node[el, below, pos=0.3] {$s_3$} (4)
[->] (2) edge [color=blue] node[el, above, pos=0.6] {$s_3$} (5);
\node[draw=none, below of =1, xshift=-0.5cm, yshift=-1.9cm] {$\overrightarrow{R}(G_3)$};
\end{scope}

\begin{scope}[xshift=8.5cm, yshift=-0.5cm]
\node [fill, label=left:$0010$] (1) {};
\node [fill, label=left:$1010$] (2) [above of=1] {};
\node [fill, label=left:$1110$] (3) [above of=2] {};
\node [fill, label=left:$0000$] (4) [below of=1,  xshift=-1cm] {};
\node [fill, label=left:$1000$] (5) [below of=2,  xshift=-1cm] {};
\node [fill, label=right:$0011$] (6) [above of=1,  xshift=1cm] {};
\node [fill, label=right:$1011$] (7) [above of=2,  xshift=1cm] {};
\node [fill, label=right:$1111$] (8) [above of=3,  xshift=1cm] {};
\path[every node/.style={font=\sffamily\small}, inner sep=0.5pt, line width=1.pt]
[<-] (1) edge [color=blue] node[el, right, pos=0.6] {$s_1$}  (2)
[<-] (4) edge [color=blue] node[el, left, pos=0.5] {$s_1$} (5)
[<-] (2) edge [color=blue] node[el, left, pos=0.5] {$s_2$} (3)
[->] (1) edge [color=blue] node[el, below, pos=0.3] {$s_3$} (4)
[->] (2) edge [color=blue] node[el, above, pos=0.6] {$s_3$} (5)
[<-] (1) edge [color=blue] node[el, below, pos=0.7] {$s_4$} (6)
[<-] (2) edge [color=blue] node[el, above, pos=0.5] {$s_4$} (7)
[<-] (3) edge [color=blue] node[el, above, pos=0.5] {$s_4$} (8)
[<-] (6) edge [color=blue] node[el, right, pos=0.5] {$s_1$} (7)
[<-] (7) edge [color=blue] node[el, right, pos=0.5] {$s_2$} (8);
\node[draw=none, below of =1, yshift=-1.2cm]{$\overrightarrow{R}(G_4)$};
\end{scope}

\begin{scope}[xshift=12cm]
\node [fill, label=left:$00101$] (1a) {};
\node [fill, label=left:$10101$] (2a) [above of=1a] {};
\node [fill, label=left:$11101$] (3a) [above of=2a] {};
\node [fill, label=left:$00001$] (4a) [below of=1a,  xshift=-1cm, yshift=0.3cm] {};
\node [fill, label=left:$10001$] (5a) [below of=2a,  xshift=-1cm, yshift=0.3cm] {};
\node [fill, label=right:$00111$] (6a) [above of=1a, xshift=1cm, yshift=-0.3cm] {};
\node [fill, label=right:$10111$] (7a) [above of=2a, xshift=1cm, yshift=-0.3cm] {};
\node [fill, label=right:$11111$] (8a) [above of=3a, xshift=1cm, yshift=-0.3cm] {};
\path[every node/.style={font=\sffamily\small}, inner sep=0.5pt, line width=1.pt]
[<-] (1a) edge [color=blue] node[el, right, pos=0.5] { $s_1$} (2a)
[<-] (4a) edge [color=blue] node[el, left, pos=0.5] { $s_1$} (5a)
[<-] (2a) edge [color=blue] node[el, left, pos=0.5] { $s_2$} (3a)
[->] (1a) edge [color=blue] node[el, below, pos=0.5] { $s_3$} (4a)
[->] (2a) edge [color=blue] node[el, above, pos=0.7] { $s_3$} (5a)
[<-] (1a) edge [color=blue] node[el, right, pos=0.1] { $s_4$} (6a)
[<-] (2a) edge [color=blue] node[el, right, pos=0.3] { $s_4$} (7a)
[<-] (3a) edge [color=blue] node[el, right, pos=0.3] { $s_4$} (8a)
[<-] (6a) edge [color=blue] node[el, right, pos=0.5] { $s_1$} (7a)
[<-] (7a) edge [color=blue] node[el, right, pos=0.5] { $s_2$} (8a);
\end{scope}

\begin{scope}[xshift=13.5cm]
\node [fill, label=right:$00100$] [yshift=-1.25cm](1b) {};
\node [fill, label=right:$10100$] (2b) [above of=1b] {};
\node [fill, label=right:$00000$] (4b) [below of=1b, xshift=-1cm, yshift=0.3cm] {};
\node [fill, label=right:$10000$] (5b) [below of=2b, xshift=-1cm, yshift=0.3cm] {};
\node [fill, label=right:$00110$] (6b) [above of=1b, xshift=1cm, yshift=-0.3cm] {};
\node [fill, label=right:$10110$] (7b) [above of=2b, xshift=1cm, yshift=-0.3cm] {};
\path[every node/.style={font=\sffamily\small}, inner sep=0.5pt, line width=1.pt]
[<-] (1b) edge [color=blue] node[el, right, pos=0.5] { $s_1$} (2b)
[<-] (4b) edge [color=blue] node[el, right, pos=0.5] { $s_1$} (5b)
[->] (1b) edge [color=blue] node[el, right, pos=0.5] { $s_3$} (4b)
[->] (2b) edge [color=blue] node[el, right, pos=0.5] { $s_3$} (5b)
[<-] (1b) edge [color=blue] node[el, right, pos=0.6] { $s_4$} (6b)
[<-] (2b) edge [color=blue] node[el, right, pos=0.6] { $s_4$} (7b)
[<-] (6b) edge [color=blue] node[el, right, pos=0.5] { $s_1$} (7b)
[<-] (1b) edge [color=blue] node[el, right, pos=0.7] { $s_5$} (1a)
[<-] (2b) edge [color=blue] node[el, right, pos=0.4] { $s_5$} (2a)
[<-] (4b) edge [color=blue] node[el, left, pos=0.4] { $s_5$} (4a)
[<-] (5b) edge [color=blue] node[el, right, pos=0.4] { $s_5$} (5a)
[<-] (6b) edge [color=blue] node[el, right, pos=0.7] { $s_5$} (6a)
[<-] (7b) edge [color=blue] node[el, right, pos=0.5] { $s_5$} (7a);
\node[draw=none, below of =1b,  xshift=-0.5cm, yshift=-0.45cm] {$\overrightarrow{R}(G)$};
\end{scope}
\end{tikzpicture}
\end{center}
\caption{\label{A-FDL} An example for implementing Algorithm \ref{A:FDL}.}
\end{figure}

See  Figure \ref{A-FDL} for the illustration of Algorithm \ref{A:FDL}. 
Note that graph $G$ in  Figure \ref{A-FDL} is a peripherally 2-colorable graph.
Peripherally 2-colorable graphs introduced in \cite{BCTZ23+} form a subclass of plane elementary bipartite graphs whose infinite face is forcing.  
For the same graph $G$, the binary coding on the vertex set of $R(G)$ as a finite distributive lattice in Figure \ref{A-FDL}  generated by Algorithm \ref{A:FDL}
 is different from the binary coding on the vertex set of $R(G)$ as a daisy cube displayed in Figure 2 in \cite{BCTZ23+}.

The following theorem is a generalization of Theorem 5.2 in \cite{ZLS08} 
for catacondensed benzenoid graphs.

\begin{theorem}\label{T:DistributiveLatticeR(G)-ThetaClass-RF} 
Let $G$ be a plane elementary bipartite graph whose infinite face is forcing.  
Assume that $G$ has a $\mathrm{RFD}(G_1, G_2, \ldots, G_n)$   
associated with a sequence of finite faces $s_1, s_2, \ldots, s_n$ for some positive integer $n$
such that $s_i$ is a reducible face of $G_i$ for all $2 \le i \le n$. 
Let $\mathcal{M}(G)$ be the set of all perfect matchings of $G$.
Then Algorithm \ref{A:FDL} determines a binary coding of length $n$ on $\mathcal{M}(G)$
such that the binary string for any $M \in \mathcal{M}(G)$ is $\phi_{M}=\phi_{M}(s_1)\phi_{M}(s_2)\ldots \phi_{M}(s_{n})$.  
\end{theorem} 
\proof   Prove by induction on $n$. If $n=1$, then it is trivial since $G$ is an even cycle 
with exactly two perfect matchings  which have binary codes $0$ and $1$, respectively. 
Let  $n \ge 2$. Since $G$ is elementary, by \cite{LZ03}, $\mathbb{M}(G)=(\mathcal{M}(G), \le_L)$ 
is a finite distributive lattice whose Hasse diagram is isomorphic to $\overrightarrow{R}(G)$.
Let $M \in \mathcal{M}(G)$. 
Since the infinite face of $G$ is forcing, by Theorem \ref{T:IsometricDimension}, 
$\phi: M \rightarrow \phi_M$ defines a binary coding 
$\phi_M=\phi_M(s_1)\phi_M(s_2) \ldots \phi_M(s_n)$ of length $n$ on $\mathcal{M}(G)$.
It follows that $\phi_M(s_i)=0$ or $1$ for all $1 \le i \le n$.

Note that $s_n$ is a reducible face of $G$. Let  $P_n$ be the common periphery of $s_n$ and $G$.
Let $H$ be the subgraph of $G$ obtained by removing all internal vertices and edges of  $P_n$.
Recall that the infinite face of $G$ is forcing. Then by Lemma \ref{L:Basic}, $R(G)$ can be constructed from 
the one-edge graph by a sequence of peripheral convex expansions with respect to a $\mathrm{RFD}(G_1, G_2, \ldots, G_n)$
associated with a sequence of finite faces $s_1, s_2, \ldots, s_n$
such that $s_i$ is a reducible face of $G_i$ for all $2 \le i \le n$.
Note that $H=G_{n-1}$.
Then $R(H)$  can be constructed from 
the one-edge graph by a sequence of peripheral convex expansions with respect to a $\mathrm{RFD}(G_1, G_2, \ldots, G_{n-1})$
associated with a sequence of finite faces $s_1, s_2, \ldots, s_{n-1}$
such that $s_i$ is a reducible face of $G_i$ for all $2 \le i \le n-1$.
By Lemma \ref{L:Basic}, the infinite face of $H$ is forcing.
Let $\mathcal{M}(H)$ be the set of all perfect matchings of $H$.
Note that $R(G)$ can be obtained from $R(H)$ by a peripheral convex expansion with respect to 
a reducible face $s_n$ of $G$. Then by Theorem \ref{T:MedianR(G)}, $\mathcal{M}(G; P_n^{+}) = \mathcal{M}(G; P_n^{+}, \partial s_n)$,
and $\mathcal{M}(G)=\mathcal{M}(G; P_n^{-}) \cup \mathcal{M}(G; P_n^{+}, \partial s_n)$.
 Let $M \in \mathcal{M}(G)$. Then there are two possible cases.
\vskip 0.1in

Case 1. $M \in \mathcal{M}(G; P_n^{-})$. Let $M^H$ be the restriction of $M$ on $H$.
It is clear that $M^H$ is a perfect matching of $H$.
Since the infinite face of $H$ is forcing, by induction hypothesis,  Algorithm \ref{A:FDL} determines a binary string $x' \in \mathbf{L_{n-1}}$ 
for $M^H \in \mathcal{M}(H)$ such that $x'=\phi_{M^H}=\phi_{M^H}(s_1) \ldots \phi_{M^H}(s_{n-1})$. 

By Lemma \ref{L:Claim1},  $\phi_M(s_i)=\phi_{M^H}(s_i)$ for all $1 \le i \le n-1$.
Moreover, if $P_n$  has two end vertices colored white and black respectively 
along the clockwise orientation of the periphery of $G$, then $\phi_M(s_n)=0$,  and so $\phi_M=x'0$; 
Otherwise, $\phi_M(s_n)=1$, and so $\phi_M=x'1$.

\vskip 0.1in
Case 2. $M \in  \mathcal{M}(G; P_n^{+}, \partial s_n)$. 
Then $M \oplus \partial s_n \in   \mathcal{M}(G; P_n^{-}, \partial s) \subseteq  \mathcal{M}(G; P_n^{-})$.
Let $(M \oplus \partial s_n)^H$ be the restriction of $M \oplus \partial s$ on $H$.
It is clear that $(M \oplus \partial s_n)^H$ is a perfect matching  of $H$. 
By induction hypothesis,  Algorithm \ref{A:FDL} determines a binary string 
$x'' \in \mathbf{L_{n-1}}$ for $(M \oplus \partial s_n)^H \in  \mathcal{M}(H)$ such that 
$x''=\phi_{(M \oplus \partial s_n)^H}
=\phi_{(M \oplus \partial s_n)^H}(s_1) \ldots \phi_{(M \oplus \partial s_n)^H}(s_{n-1})$. 
 
 By Lemma \ref{L:Claim2},  $\phi_M(s_i)=\phi_{(M \oplus \partial s_n)^H}(s_i)$ for all $1 \le i \le n-1$. 
Let $s_{\alpha(n,1)}, \ldots, s_{\alpha(n,t_n)}$ be the finite faces of $G$ having common edges
with $s_n$, where $\alpha(n,j)$ is an integer in $\{1, 2, \ldots, n-1\}$ for all $j \in \{1, 2, \ldots, t_n\}$.
Then by Lemma \ref{L:Claim2},  if $P_n$  has two end vertices colored white and black respectively 
along the clockwise orientation of the periphery of $G$, 
then $\phi_M(s_n)=\phi_{M\oplus \partial s_n} (s_{\alpha(n,j)})=1$ for all $1 \le j \le n_t$, 
and  $\phi_{M}=x''1$; Otherwise,
$\phi_M(s_n)=\phi_{M\oplus \partial s_n} (s_{\alpha(n,j)})=0$ for all $1 \le j \le n_t$,  
and $\phi_{M}=x''0$.\qed\\

\vskip 0.1in
\noindent {\bf Assumption 1.}
Let $G$ be a plane elementary bipartite graph whose infinite face is forcing.  
Assume that $G$ has a $\mathrm{RFD}(G_1, G_2, \ldots, G_n)$   
associated with a sequence of finite faces $s_1, s_2, \ldots, s_n$ for some positive integer $n$
such that $s_i$ is a reducible face of $G_i$ for all $2 \le i \le n$.
Let $s_{\alpha(n,1)}, \ldots, s_{\alpha(n,t_n)}$ be the finite faces of $G$ having common edges
with $s_n$, where $\alpha(n,j)$ is an integer in $\{1, 2, \ldots, n-1\}$ for all $j \in \{1, 2, \ldots, t_n\}$.
Let $P_n$ be the common periphery of $s_n$ and $G$.
Let  $H$ be the subgraph of $G$ obtained by removing all internal vertices and edges from $P_n$.

\begin{lemma}\label{L:Claim1} 
Let Assumption 1 hold true.
Let $M \in \mathcal{M}(G; P_n^{-})$ and $M^H$ be the restriction of $M$ on $H$.
Then $M^H$ is a perfect matching of $H$, and $\phi_M(s_i)=\phi_{M^H}(s_i)$ for all $1 \le i \le n-1$.
Moreover, $\phi_M(s_n)=0$ if $P_n$  has two end vertices colored white and black respectively 
along the clockwise orientation of the periphery of $G$; and $\phi_M(s_n)=1$ otherwise.
\end{lemma} 

\proof Note that  $s_n$ is a reducible face of $G$.
If $s_{\alpha(n,1)}, \ldots, s_{\alpha(n,t_n)}$ are the finite faces of $G$ having common edges
with $s_n$, where $\alpha(n,j)$ is an integer in $\{1, 2, \ldots, n-1\}$ for all $j \in \{1, 2, \ldots, t_n\}$,
then the facial cycle $\partial s_n=P_n \cup J_{\alpha(n)}$, where $P_n$ is an odd length handle of $G$ on $\partial G$,
and $J_{\alpha(n)}$ is an odd length path of $G$  on $\partial H$. 
Moreover, $P_n$ and $J_{\alpha(n)}$ have two end vertices in common.

We have seen that both $G$ and $H$ are elementary. 
Let $\mathcal{M}(G)$ (respectively, $\mathcal{M}(H)$) be the set of all perfect matchings of $G$ (respectively, $H$).
Then by \cite{LZ03}, 
$\mathbb{M}(G)=(\mathcal{M}(G), \le_L)$ and $\mathbb{M}(H)=(\mathcal{M}(H), \le_L)$ 
are finite distributive lattices.
Let $M_{\hat{0}}$ be  the minimum element of $\mathbb{M}(G)$, and
$M_{[0]}$ be the minimum element of  $\mathbb{M}(H)$.
Note that $M_{[0]}$ can be extended to a unique perfect matching $\overline{M_{[0]}}$ in $\mathcal{M}(G; P_n^{-})$. 
Let $M \in \mathcal{M}(G; P_n^{-})$ and $M^H$ be the restriction of $M$ on $H$.
It is easily seen that $M^H$ is a perfect matching of $H$.

We first show that  $\phi_M(s_i)=\phi_{M^H}(s_i)$ for all $1 \le i \le n-1$
based on whether $\overline{M_{[0]}}$ is equal to $M_{\hat{0}}$ or not.

Case 1.  $\overline{M_{[0]}}=M_{\hat{0}}$.
Then $M_{\hat{0}} \in \mathcal{M}(G; P_n^{-})$ and $M_{[0]}=M^H_{\hat{0}}$ which is the restriction of $M_{\hat{0}}$ on $H$.
It follows that  $M \oplus M_{\hat{0}}$ does not contain any edges on $P_n$ 
since both $M$ and $M_{\hat{0}}$ are contained in $\mathcal{M}(G; P_n^{-})$.
Hence, $M \oplus M_{\hat{0}}$ and $M^H \oplus M^H_{\hat{0}}$ are the same. 
By Definition \ref{D:Functions}, 
for all $1 \le i \le n-1$, 
$\phi_M(s_i)$  is the number of cycles in $M \oplus M_{\hat{0}}$ with $s_i$ in their interiors;
and $\phi_{M^H}(s_i)$  is the number of cycles in $M^H \oplus M^H_{\hat{0}}$  with $s_i$ in their interiors.
It follows that $\phi_M(s_i)=\phi_{M^H} (s_i)$ for all $1 \le i \le n-1$.

Case 2. $\overline{M_{[0]}} \neq M_{\hat{0}}$. 
Then $M_{[0]}$ is not the restriction of $M_{\hat{0}}$ on $H$.
Recall that  $M_{\hat{0}}$ is the minimum of $\mathbb{M}(G)$,
and $M_{[0]}$ is the minimum of $\mathbb{M}(H)$.
Then $M_{\hat{0}}$ the unique perfect matching $G$ 
such that $G$ has no proper  $M_{\hat{0}}$-alternating cycles,
and $M_{[0]}$ 
is the unique perfect matching $H$ such that $H$ has no proper  $M_{[0]}$-alternating cycles.
By \cite{TV12},  $\partial G$ is an even cycle that is $M_{\hat{0}}$-alternating,
and  $\partial H$ is an even cycle that is  $M_{[0]}$-alternating.
Recall that $P_n$ and $J_{\alpha(n)}$ have two end vertices in common,
where $P_n$ is an odd length handle of $G$ on $\partial G$, and $J_{\alpha(n)}$ is an odd length path of $G$ on $\partial H$.
Then $M_{\hat{0}} \in \mathcal{M}(G; P_n^{+})$ and $M_{[0]} \in \mathcal{M}(H; J_{\alpha(n)}^{+})$.

Note that $\mathcal{M} (G; P_n^{-})$ and $\mathcal{M} (G; P_n^{+})$ form a partition of $\mathcal{M}(G)$
since $P_n$ is an odd length handle of $G$ on $\partial G$.   Then $M_{[0]} \in \mathcal{M}(H; J_{\alpha(n)}^{+})$
implies that  $\overline{M_{[0]}} \in \mathcal{M}(G; P_n^{-}) \cap \mathcal{M}(G; J_{\alpha(n)}^{+})$
since $P_n$ and $J_{\alpha(n)}$ have two end vertices in common.
We will show that $M_{\hat{0}} \in \mathcal{M}(G; P_n^{+})$  
implies that $M_{\hat{0}} \in \mathcal{M}(G; P_n^{+}) \cap \mathcal{M}(G; J_{\alpha(n)}^{-})$.
Recall that the infinite face of $G$ is forcing. Then by Lemma \ref{L:Basic}, $R(G)$ can be constructed from 
the one-edge graph by a sequence of peripheral convex expansions with respect to a $\mathrm{RFD}(G_1, G_2, \ldots, G_n)$
associated with a sequence of finite faces $s_1, s_2, \ldots, s_n$
such that $s_i$ is a reducible face of $G_i$ for all $2 \le i \le n$.
In particular, $R(G)$ can be obtained from $R(H)$ by a peripheral convex expansion with respect to the reducible face $s_n$ of $G$.
By Theorem \ref{T:MedianR(G)}, $\mathcal{M}(G; P_n^{+})= \mathcal{M}(G; P_n^{+}, \partial s_n)$.
Recall that $\partial s_n=P_n \cup J_{\alpha(n)}$, where $J_{\alpha(n)}$ is an odd length path of $G$, 
and $P_n$ and  $J_{\alpha(n)}$ have two end vertices in common.
Then $M_{\hat{0}} \in \mathcal{M}(G; P_n^{+})= \mathcal{M}(G; P_n^{+}, \partial s_n)$ implies that $J_{\alpha(n)}$ 
is an  $M_{\hat{0}}$-alternating path without end edges contained in $M_{\hat{0}}$.
Hence, $M_{\hat{0}} \in \mathcal{M}(G; P_n^{+}) \cap \mathcal{M}(G; J_{\alpha(n)}^{-})$.

Now, we have $\overline{M_{[0]}} \in \mathcal{M}(G; P_n^{-}) \cap \mathcal{M}(G; J_{\alpha(n)}^{+})$
and $M_{\hat{0}} \in \mathcal{M}(G; P_n^{+}) \cap \mathcal{M}(G; J_{\alpha(n)}^{-})$.
So, $P_n \cup J_{\alpha(n)}$ is contained in $\overline{M_{[0]}} \oplus M_{\hat{0}}$.
It is well known \cite{LP86} that $\overline{M_{[0]}} \oplus M_{\hat{0}}$ is a set of vertex disjoint 
$(\overline{M_{[0]}}, M_{\hat{0}})$-alternating cycles.  
Note that the restriction of $\overline{M_{[0]}}$ on $H$ is $M_{[0]}$, which is the minimum of $\mathbb{M}(H)$.
Then $H$ does not contain any proper $M_{[0]}$-alternating cycles,
and so $H$ does not contain any proper $\overline{M_{[0]}}$-alternating cycles either.
Note that $G$ does contain any proper $M_{\hat{0}}$-alternating cycles since $M_{\hat{0}}$ is the minimum of $\mathbb{M}(G)$,
and so $H$ doe not contain any proper $M_{\hat{0}}$-alternating cycles either.
It follows that there are no $(\overline{M_{[0]}}, M_{\hat{0}})$-alternating cycles in $H$.
Hence, $\overline{M_{[0]}} \oplus M_{\hat{0}}$ contains a unique $(\overline{M_{[0]}}, M_{\hat{0}})$-alternating cycle
which is the facial cycle $\partial s_n=P_n \cup J_{\alpha(n)}$, that is,
$\overline{M_{[0]}} \oplus M_{\hat{0}}=\partial s_n$.
Then $\overline{M_{[0]}}$ covers $M_{\hat{0}}$ since $M_{\hat{0}}$ is the minimum of $\mathbb{M}(G)$. 
By Lemma \ref{L:phi-psi}, 
 $\phi_{\overline{M_{[0]}}}(s_i)-\phi_{M_{\hat{0}}}(s_i)=0$ for all $1 \le i \le n-1$.
Moreover, $\phi_{M_{\hat{0}}}(s_i)=0$  since $M_{\hat{0}}$ is the minimum of $\mathbb{M}(G)$. 
Then $\phi_{\overline{M_{[0]}}}(s_i)=\phi_{M_{\hat{0}}}(s_i)=0$ for all $1 \le i \le n-1$.

By Lemma \ref{L:phi-psi}, for any perfect matching $M$ of $G$, we have 
$\phi_M - \phi_{M_{\hat{0}}} =\psi_{MM_{\hat{0}}}$, and $\phi_M-\phi_{\overline{M_{[0]}}}=\psi_{M\overline{M_{[0]}}}$.
Use the above fact that $\phi_{\overline{M_{[0]}}}(s_i)=\phi_{M_{\hat{0}}}(s_i)=0$   for all $1 \le i \le n-1$,
we have \[\phi_M(s_i)=\psi_{MM_{\hat{0}}}(s_i)=\psi_{M\overline{M_{[0]}}}(s_i) \  \mathrm{ for \  all } \ 1 \le i \le n-1. \ \ \ \ (1) \]

Recall that $M \in \mathcal{M}(G; P_n^{-})$ and $\overline{M_{[0]}} \in \mathcal{M}(G; P_n^{-})$.
Then $M \oplus \overline{M_{[0]}}$ does not contain any edges on $P_n$. 
Note that $M^H$ is the restriction of $M$ on $H$,
and $M_{[0]}$ is the restriction of $\overline{M_{[0]}}$ on $H$.
Then $M \oplus \overline{M_{[0]}}$ and $M^H \oplus M_{[0]}$ are the same.
By Definition \ref{D:Functions},   we have
\[\psi_{M\overline{M_{[0]}}}(s_i)=\psi_{M^HM_{[0]}}(s_i)  \  \mathrm{ for \  all } \ 1 \le i \le n-1. \ \ \ \ (2) \]

By Lemma \ref{L:phi-psi}, 
 $\psi_{M^HM_{[0]}}(s_i)=\phi_{M^H}(s_i)- \phi_{M_{[0]}}(s_i)$  for all $1 \le i \le n-1$.
Moreover, $\phi_{M_{[0]}}(s_i)=0$ since $M_{[0]}$ is the minimum of $\mathbb{M}(H)$. 
It follows that
\[\psi_{M^HM_{[0]}}(s_i)=\phi_{M^H}(s_i)  \  \mathrm{ for \  all } \ 1 \le i \le n-1. \ \ \ \ (3) \]

By Equations (1) - (3), we have 
\[\phi_M(s_i)=\psi_{MM_{\hat{0}}}(s_i)=\psi_{M\overline{M_{[0]}}}(s_i)
=\psi_{M^HM_{[0]}}(s_i)=\phi_{M^H}(s_i)  \  \mathrm{ for \  all } \ 1 \le i \le n-1. \]

Next, we determine $\phi_M(s_n)$ based on the colors of  two end vertices of $P_n$
along the clockwise orientation of the periphery of $G$.
Recall that $P_n$ is on the periphery of $G$ which is improper $M_{\hat{0}}$-alternating.

If $P_n$  has two end vertices colored white and black respectively 
along the clockwise orientation of the periphery of $G$,
then $M_{\hat{0}}  \in \mathcal{M}(G; P_n^{-})$.
By our assumption that $M \in \mathcal{M}(G; P_n^{-})$,
it follows that $M \oplus M_{\hat{0}}$ does not contain any edges on $P_n$.
It is well known \cite{LP86} that $M \oplus M_{\hat{0}}$ is a set of vertex disjoint $(M, M_{\hat{0}})$-alternating cycles of $G$.
So, $s_n$ is not contained in the interior region of any cycle in $M \oplus M_{\hat{0}}$.
This implies that  $\phi_M(s_n)=0$.  

If $P_n$  has two end vertices colored black and white respectively 
along the clockwise orientation of the periphery of $G$,
then $M_{\hat{0}}  \in \mathcal{M}(G; P_n^{+})$.
By our assumption that $M \in \mathcal{M}(G; P_n^{-})$,
it follows that $P_n$ is contained in $M \oplus M_{\hat{0}}$.
Hence, $s_n$ is contained in the interior region of a cycle in $M \oplus M_{\hat{0}}$,
and so $\phi_M(s_n) \ge 1$ by Definition \ref{D:Functions}.
By Theorem \ref{T:IsometricDimension}, 
$\phi: M \rightarrow \phi_M$ defines a binary coding 
$\phi_M=\phi_M(s_1)\phi_M(s_2) \ldots \phi_M(s_n)$ of length $n$ on $\mathcal{M}(G)$.
It follows that  $\phi_M(s_n)=1$.
  \qed\\

\begin{lemma}\label{L:Claim2} 
Let Assumption 1 hold true.
Let $M \in \mathcal{M}(G; P_n^{+}, \partial s_n)$, and $(M \oplus \partial s_n)^H$ be the restriction of $M \oplus \partial s$ on $H$.
Then $(M \oplus \partial s_n)^H$ is a perfect matching of $H$ and $\phi_M(s_i)=\phi_{(M \oplus \partial s_n)^H}(s_i)$ for all $1 \le i \le n-1$. 
Moreover,  if $P_n$  has two end vertices colored white and black respectively along the clockwise orientation of the periphery of $G$,
then $\phi_{M}(s_n)=\phi_{M\oplus \partial s_n} (s_{\alpha(n,j)})=1$ for all $1 \le j \le n_t$; 
Otherwise, $\phi_{M}(s_n)=\phi_{M\oplus \partial s_n} (s_{\alpha(n,j)})=0$ for all $1 \le j \le n_t$.
\end{lemma}
\proof Let  $M \in \mathcal{M}(G; P_n^{+}, \partial s_n)$. 
Then $M \oplus \partial s_n \in \mathcal{M}(G; P_n^{-}, \partial s_n)  \subseteq \mathcal{M}(G; P^{-})$.
Note that $M \oplus (M \oplus \partial s_n) = \partial s_n$.
Then there is an edge between $M$ and $M \oplus \partial s_n$ in $R(G)$ with the face-label $s_n$.
So, either $M$ covers $M \oplus \partial s_n$, or vice versa.
By Lemma \ref{L:phi-psi}, $\phi_M(s_n)-\phi_{M \oplus \partial s_n}(s_n)= \pm 1$, 
and $\phi_{M}(s_i)=\phi_{M \oplus \partial s_n}(s_i)$ for $1 \le i \le n-1$.
Let $(M \oplus \partial s_n)^H$ be the restriction of $M \oplus \partial s_n$ on $H$. 
Since $M \oplus \partial s_n \in \mathcal{M}(G; P_n^{-})$,
by Lemma \ref{L:Claim1},  $(M \oplus \partial s_n)^H$ is a perfect matching of $H$,
and $\phi_{M \oplus \partial s_n}(s_i)=\phi_{(M \oplus \partial s_n)^H} (s_i)$ for all $1 \le i \le n-1$. 
It follows that $\phi_M(s_i)=\phi_{(M \oplus \partial s_n)^H}(s_i)$ for all $1 \le i \le n-1$. 

Next, we determine $\phi_M(s_n)$ based on the colors of  two end vertices of $P_n$
along the clockwise orientation of the periphery of $G$.
Assume that $P_n$  has two end vertices colored white and black respectively 
along the clockwise orientation of the periphery of $G$. 
Recall that $M \in \mathcal{M}(G; P_n^{+}, \partial s_n)$. Then $\partial s_n$ is proper $M$-alternating,
Hence, $M \in \mathcal{M}(G; P_n^{+}, \partial s_n)$ covers $M \oplus \partial s_n \in \mathcal{M}(G; P_n^{-}, \partial s_n)$.
By Lemma \ref{L:phi-psi}, $\phi_M(s_n)-\phi_{M \oplus \partial s_n}(s_n)= 1$.
It follows that  $\phi_M(s_n)=1$ since   $\phi_M$ is a binary string by Theorem \ref{T:IsometricDimension}.

By the fact that $P_n$ and $J_{\alpha(n)}$
have two end vertices in common, and the assumption on the colors of two end vertices of $P_n$,
it follows that $J_{\alpha(n)}$ has two end vertices colored white and black respectively 
along the clockwise orientation of the periphery of $H$.
Let $M_{[0]}$ be the minimum element of the finite distributive lattice $\mathbb{M}(H)=(\mathcal{M}(H), \le_L)$.
Note that $\partial H$ 
is an improper $M_{[0]}$-alternating cycle by  \cite{TV12}, 
and $J_{\alpha(n)}$ is an odd length path of $G$  on $\partial H$.
Then $M_{[0]} \in  \mathcal{M} (H; J_{\alpha(n)}^{-})$.

Recall that $P_n$ is an odd length handle of $G$ on $\partial G$, and $J_{\alpha(n)}$ is an odd length path  of $G$ on $\partial H$
such that $\partial s_n=P_n \cup J_{\alpha(n)}$ where $P_n$ and $J_{\alpha(n)}$ have two end vertices in common.
Then $\langle  \mathcal{M}(G; P_n^{-}, \partial s_n)  \rangle=\langle \mathcal{M}(G; J_{\alpha(n)}^{+}) \rangle$.
Moreover, $\langle \mathcal{M}(G; J_{\alpha(n)}^{+}) \rangle  \cong \langle \mathcal{M}(H; J_{\alpha(n)}^{+}) \rangle$
since $H$ is the subgraph of $G$ obtained by removing all internal vertices and edges from $P_n$.
By the fact that $M \oplus \partial s_n \in \mathcal{M}(G; P_n^{-}, \partial s_n)$, we have
that $(M \oplus \partial s_n)^H \in \mathcal{M} (H; J_{\alpha(n)}^{+})$.

Now, $M_{[0]} \in  \mathcal{M} (H; J_{\alpha(n)}^{-})$ and  $(M \oplus \partial s_n)^H \in \mathcal{M} (H; J_{\alpha(n)}^{+})$.
So, $(M \oplus \partial s_n)^H \oplus M_{[0]}$ contains the odd length path $J_{\alpha(n)}$.
Recall that $s_{\alpha(n,1)}, \ldots, s_{\alpha(n,t_n)}$ are the finite faces of $G$ having common edges
with $s_n$ where $1 \le \alpha(n,j) \le n-1$ for all $j \in \{1, 2, \ldots, t_n\}$.
Let $J_{\alpha(n,j)}$ be  the common periphery of $s_{\alpha(n,j)}$ and $s_n$ for all $1 \le j \le t_n$.
Then $J_{\alpha(n)}=\cup_{j=1}^{n_t} J_{\alpha(n,j)}$, where $J_{\alpha(n,j)}$ is a handle of $G$ on $\partial H$  for all $1 \le j \le t_n$.
Hence,   $(M \oplus \partial s_n)^H \oplus M_{[0]}$ contains each handle $J_{\alpha(n,j)}$ for all $1 \le j \le t_n$.
It is well known \cite{LP86} that $(M \oplus \partial s_n)^H \oplus M_{[0]}$ is a set of vertex disjoint $((M \oplus \partial s_n)^H, M_{[0]})$-alternating cycles.
So, $s_{\alpha(n,j)}$ is contained in a cycle of $(M \oplus \partial s_n)^H \oplus M_{[0]}$ for all $1 \le j \le t_n$.
Since $\phi_{(M \oplus \partial s_n)^H}$ is a binary string
by Theorem \ref{T:IsometricDimension}, it follows that $\phi_{(M \oplus \partial s_n)^H}(s_{\alpha(n,j)}) =1$  for all $1 \le j \le t_n$.
By Lemma \ref{L:Claim1}, $\phi_{M \oplus \partial s_n}(s_{\alpha(n,j)})=\phi_{(M \oplus \partial s_n)^H}(s_{\alpha(n,j)})$
since $M \oplus \partial s_n \in \mathcal{M}(G; P_n^{-})$ and $1 \le \alpha(n,j)  \le n-1$.
It follows that  $\phi_{M \oplus \partial s_n}(s_{\alpha(n,j)}) =1$ for all $1 \le j \le t_n$.

Therefore, if $P_n$  has two end vertices colored white and black respectively 
along the clockwise orientation of the periphery of $G$, then  
$\phi_{M}(s_n)=\phi_{M\oplus \partial s_n} (s_{\alpha(n,j)})=1$  for all $1 \le j \le t_n$.

Similarly, we can show that if $P_n$  has two end vertices colored black and white respectively 
along the clockwise orientation of the periphery of $G$, then 
$\phi_M(s_n)=\phi_{M\oplus \partial s_n} (s_{\alpha(n,j)})=0$  for all $1 \le j \le t_n$.
\qed\\

\subsection{When $G$ is a plane bipartite graph}  

We are ready to provide general versions of our main results.
Let $\mathbb{P}=(P, \le)$ and $\mathbb{Q}=(Q, \le)$ be two posets. The direct product of $\mathbb{P}$ and $\mathbb{Q}$,
denoted by $\mathbb{P} \times \mathbb{Q}$, consists of the ordered pairs $(p,q)$, where $p \in P$ and $q \in Q$ such that
$(p_1,q_1) \le (p_2,q_2)$ if $p_1 \le p_2$ in $\mathbb{P}$ and $q_1 \le q_2$ in $\mathbb{Q}$.
It is known \cite{Z10} that for a plane bipartite graph with a perfect matching, 
$\mathbb{M}(G)=(\mathcal{M}(G), \le_L)$ is a finite distributive lattice whose Hasse diagram 
is isomorphic to $\overrightarrow{R}(G) $ if and only if $G$ is weakly elementary.
Moreover, if $E_1, E_2, \ldots, E_k$ are the elementary components of a plane weakly elementary bipartite graph $G$, then
$\mathbb{M}(G) \cong \mathbb{M}(E_1) \times \ldots \times \mathbb{M}(E_k)$, where $\mathbb{M}(E_i)=(\mathcal{M}(E_i), \le_L)$
for all $1 \le i \le k$.

A Cartesian product of two graphs $G$ and $H$ is a graph $G \Box H$ with vertex set $V(G) \times V(H)$
such that two vertices $(g_1,h_1)$ and $(g_2,h_2)$ are adjacent in $G \Box H$ if either $g_1g_2 \in E(G)$ and $h_1=h_2$,
or  $h_1h_2 \in E(H)$ and $g_1=g_2$.  
It was shown \cite{O08} that a Cartesian product $G \Box H$ of two finite partial cubes $G$ and $H$ has the isometric dimension
$\mathrm{idim}(G \Box H)=\mathrm{idim}(G) + \mathrm{idim}(H)$.
If $G$ is a plane weakly elementary bipartite graph
whose elementary components  are $E_1, \ldots, E_k$,
then  it is well known that $R(G)=R(E_1) \Box R(E_2) \Box \cdots \Box R(E_k)$,
which can be seen easily  by the definitions of a resonance graph and a Cartesian product.

\begin{theorem}\label{T:G-IsometricDimension}
Let $G$ is a plane bipartite graph whose resonance graph $R(G)$ is connected. Let $d$ be
the number of finite faces of $G$ without forbidden edges on their peripheries.
Then $idim(R(G)) \ge d$ and the equality holds if and only if  
$G$ is a plane weakly elementary bipartite graph whose each elementary component
with more than two vertices holds the property that the infinite face is forcing.
\end{theorem}
\proof By \cite{F03, ZZY04}, the resonance graph $R(G)$ of a plane bipartite graph $G$ is connected if and only if $G$ is weakly elementary.
It follows that $G$ is a plane weakly elementary bipartite graph.
Let $E_1, E_2, \ldots, E_k$ be the elementary components of $G$ with more than two vertices obtained by deleting all forbidden edges of $G$.
Then for each $1 \le i \le k$, $E_i$ is a plane elementary bipartite graph with $n_i$ finite faces for some positive integer $n_i$,
and  $R(G)=R(E_1) \Box R(E_2) \Box \cdots \Box R(E_k)$. 
By the definition of a plane weakly elementary bipartite graph $G$, the removal of all forbidden edges of $G$ does not resulted in any new finite faces.
If $d$ is the number of finite faces of $G$ without forbidden edges on their peripheries, then $d=\sum\limits_{i=1}^{k} n_i$.

Since $E_i$ is a plane elementary bipartite graph with $n_i$ finite faces for some positive integer $n_i$,  by Theorem \ref{T:IsometricDimension}, 
$\mathrm{idim} (R(E_i)) \ge n_i$ and the equality holds if and only if the infinite face of $E_i$ is forcing for each $1 \le i \le k$.
By  \cite{ZLS08}, any connected $R(G)$ is a median graph and can be isometrically embedded into a hypercube.
By \cite{O08}, we have  $\mathrm{idim}(R(G))=\sum\limits_{i=1}^{k} \mathrm{idim}(R(E_i)) \ge \sum\limits_{i=1}^{k} n_i=d$,
and the equality holds  if and only if the infinite face of $E_i$ is forcing for all $1 \le i \le k$.
\qed\\

\begin{theorem}\label{T:General}
Let $G$ be a plane weakly elementary bipartite graph whose 
elementary components with more than two vertices are $E_1, \ldots, E_k$.
Let $d$ be the number of finite faces of $G$ without forbidden edges on their peripheries.
Let $ \phi_{M^{[E_i]}}$ be the function defined in Definition \ref{D:Functions} for any $M^{[E_i]} \in \mathcal{M}(E_i)$, 
and $\phi^{[E_i]}: M^{[E_i]} \rightarrow \phi_{M^{[E_i]}}$ for all $1 \le i \le k$.
If the infinite face of $E_i$ is forcing for all $1 \le i \le k$,
then $[\phi]=(\phi^{[E_1]}, \ldots, \phi^{[E_k]})$ defines an isometric embedding of $R(G)$ into a hypercube $Q_d$ as a finite distributive lattice 
whose Hasse diagram is isomorphic to $\overrightarrow{R}(G)$,
where  $d=idim(R(G))$ which is equal to the summation of the numbers of  finites faces of $E_i$ for all $1 \le i \le k$.
\end{theorem}
\proof It is well known that
$R(G)=R(E_1) \Box R(E_2) \Box \cdots \Box R(E_k)$.
By Theorem \ref{T:G-IsometricDimension}, we have
$d=\mathrm{idim}(R(G))=\sum\limits_{i=1}^{k} \mathrm{idim}(R(E_i))=\sum\limits_{i=1}^{k} n_i$.

By \cite{LZ03}, $\mathbb{M}(G)=(\mathcal{M}(G), \le_L)$ is a finite distributive lattice
 whose Hasse diagram is isomorphic to  $\overrightarrow{R}(G)$,  
 $\mathbb{M}(E_i)=(\mathcal{M}(E_i), \le_L)$ is a finite distributive lattice
 whose Hasse diagram is isomorphic to $\overrightarrow{R}(E_i)$ for all $1 \le i \le k$.

If the infinite face of $E_i$ is forcing for all $1 \le i \le k$, then by lemma \ref{L:Basic}, $\phi^{[E_i]}: M^{[E_i]} \rightarrow \phi_{M^{[E_i]}}$ defines
an isometric embedding of $R(E_i)$ into a hypercube $Q_{n_i}$ as a finite distributive lattice $\mathbb{M}(E_i)=(\mathcal{M}(E_i), \le_L)$ for all $1 \le i \le k$.
By \cite{Z10},   $\mathbb{M}(G) =\mathbb{M}(E_1) \times \ldots \times \mathbb{M}(E_k)$.
It follows that $[\phi]=(\phi^{[E_1]}, \ldots, \phi^{[E_k]})$ defines an isometric embedding of $R(G)$ into a hypercube $Q_d$ as a finite distributive lattice 
whose Hasse diagram is isomorphic to $\overrightarrow{R}(G)$\qed\\

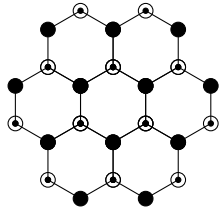
\begin{figure}[h!]
\begin{center}
\tikzset{Bullet/.style={draw=black, circle, minimum size=3pt, scale=0.5}}
\begin{tikzpicture}[scale=0.5]
  
  \draw [color=black, mark=*] plot[samples at={-150,-90,...,150,210},  variable=\x] 
  (\x:1);
    \node [Bullet] at (-150:1){};
  \node  [Bullet, fill] at (-90:1) {};
  \node [Bullet] at (-30:1) {} ;
  \node [Bullet, fill] at (30:1) {} ;
  \node  [Bullet] at (90:1) {};
  \node [Bullet, fill] at (150:1) {};
  \node [blue]{};
  
   \begin{scope}[xshift=1.732 cm]
    \draw [color=black,mark=*] plot[samples at={-150,-90,...,210},variable=\x] 
  (\x:1);
      \node [Bullet] at (-150:1){};
  \node  [Bullet, fill] at (-90:1) {};
  \node [Bullet] at (-30:1) {} ;
  \node [Bullet, fill] at (30:1) {} ;
  \node  [Bullet] at (90:1) {};
  \node [Bullet, fill] at (150:1) {};
  \node[blue]{};
  \end{scope}

 \begin{scope}[xshift=0.866 cm,  yshift=-1.5cm]
    \draw [color=black,mark=*] plot[samples at={-150,-90,...,210},variable=\x] 
  (\x:1);
      \node [Bullet] at (-150:1){};
  \node  [Bullet, fill] at (-90:1) {};
  \node [Bullet] at (-30:1) {} ;
  \node [Bullet, fill] at (30:1) {} ;
  \node  [Bullet] at (90:1) {};
  \node [Bullet, fill] at (150:1) {};
  \node[blue]{};
  \end{scope}

   \begin{scope}[xshift=-0.866 cm,  yshift=-1.5cm]
    \draw [color=black,mark=*] plot[samples at={-150,-90,...,210},variable=\x] 
  (\x:1);
      \node [Bullet] at (-150:1){};
  \node  [Bullet, fill] at (-90:1) {};
  \node [Bullet] at (-30:1) {} ;
  \node [Bullet, fill] at (30:1) {} ;
  \node  [Bullet] at (90:1) {};
  \node [Bullet, fill] at (150:1) {};
  \node[blue]{};
  \end{scope}

     \begin{scope}[xshift=2.598 cm,  yshift=-1.5cm]
    \draw [color=black,mark=*] plot[samples at={-150,-90,...,210},variable=\x] 
  (\x:1);
      \node [Bullet] at (-150:1){};
  \node  [Bullet, fill] at (-90:1) {};
  \node [Bullet] at (-30:1) {} ;
  \node [Bullet, fill] at (30:1) {} ;
  \node  [Bullet] at (90:1) {};
  \node [Bullet, fill] at (150:1) {};
  \node[blue]{};
  \end{scope}

  \begin{scope}[yshift=-3cm]
    \draw [color=black,mark=*] plot[samples at={-150,-90,...,210},variable=\x] 
  (\x:1);
      \node [Bullet] at (-150:1){};
  \node  [Bullet, fill] at (-90:1) {};
  \node [Bullet] at (-30:1) {} ;
  \node [Bullet, fill] at (30:1) {} ;
  \node  [Bullet] at (90:1) {};
  \node [Bullet, fill] at (150:1) {};
   \node [blue]{};
  \end{scope}
  
    \begin{scope}[xshift=1.732 cm, yshift=-3cm]
    \draw [color=black,mark=*] plot[samples at={-150,-90,...,210},variable=\x] 
  (\x:1);
      \node [Bullet] at (-150:1){};
  \node  [Bullet, fill] at (-90:1) {};
  \node [Bullet] at (-30:1) {} ;
  \node [Bullet, fill] at (30:1) {} ;
  \node  [Bullet] at (90:1) {};
  \node [Bullet, fill] at (150:1) {};
   \node [blue]{};
  \end{scope}  
\end{tikzpicture}
\end{center}
\caption{\label{Coronene} Coronene.}
\end{figure}

\begin{corollary}\label{C:Benzenoid}
Let $H$ be a benzenoid system whose  elementary components with more than two vertices are $H_1, \ldots, H_k$.
Let $d$ be the number of finite faces of $H$ without forbidden edges on their peripheries.
Let $ \phi_{M^{[H_i]}}$ be the function defined in Definition \ref{D:Functions} for any $M^{[H_i]} \in \mathcal{M}(H_i)$, 
and $\phi^{[H_i]}: M^{[H_i]} \rightarrow \phi_{M^{[H_i]}}$ for all $1 \le i \le k$.
If $H$ has no coronenes as its nice subgraphs,
then $[\phi]=(\phi^{[H_1]}, \ldots, \phi^{[H_k]})$ defines an isometric embedding of $R(H)$ into a hypercube $Q_d$ as a finite distributive lattice 
whose Hasse diagram is isomorphic to $\overrightarrow{R}(H)$,
where  $d=idim(R(H))$ which is equal to the summation of the numbers of  finites faces of $H_i$ for all $1 \le i \le k$.
\end{corollary}
\proof  Let $H$ be a benzenoid system.  Then $H$ is a plane weakly elementary bipartite graph \cite{ZZ00}.
If $H$ has no coronenes as its nice subgraphs, then it is trivial that
each elementary component of $H$ has no coroenes as  nice subgraphs.
On other hand, if $H$ has a coronene $T$ as its nice subgraph,
then by the definition of a nice subgraph, each perfect matching of the coronene  $T$ can be extended to a perfect matching of $H$.
Note that each edge of a coronene $T$ is contained in a perfect matching of $T$,
and so contained in a perfect matching of $H$. This implies that the coronene $T$ is contained in an elementary component of $H$.
Therefore, a benzenoid system $H$ has no coronenes as its nice subgraphs if and only if each elementary component of $H$ has no coroenes as nice subgraphs.
Let $H_1, \ldots, H_k$ be the elementary components of $H$ such that each $H_i$ has more than two vertices for all $1 \le i \le k$. 
Then each $H_i$ is an elementary benzenoid system that has no coronenes as its nice subgraphs  for all $1 \le i \le k$.
By Corollary 3.3 in \cite{BCTZ24+}, an elementary benzenoid system  $H_i$
 has no coronenes as its nice subgraphs if and only if the infinite face of $H_i$ is forcing for all $1 \le i \le k$.
 The conclusion follows by Theorem \ref{T:General}.
\qed\\

\begin{theorem}\label{T:A-General}
Let $G$ be a plane weakly elementary bipartite graph whose 
elementary components with more than two vertices are $E_1, \ldots, E_k$.
Let $d$ be the number of finite faces of $G$ without forbidden edges on their peripheries,
and $n_i$ be the number of finite faces of $E_i$ for all $1 \le i \le k$.
If the infinite face of $E_i$ is forcing for all $1 \le i \le k$, then Algorithm \ref{A:FDL} can be applied to $E_i$ 
to generate a binary coding $\mathbf{L_{n_i}^{[E_i]}}$ of length $n_i$ for all perfect matchings of $E_i$ where $1 \le i \le k$.
A binary coding $\mathbf{L_d}$ of length $d$ for all perfect matchings of $G$ can be obtained
such that any $x \in \mathbf{L_d}$ is $x=x^{[E_1]} \ldots x^{[E_k]}$ where $x^{[E_i]} \in \mathbf{L_{n_i}^{[E_i]}}$ for all $1 \le i \le k$,
inducing an isometric embedding of $R(G)$ into a hypercube $Q_d$ as a finite distributive lattice 
whose Hasse diagram is isomorphic to $\overrightarrow{R}(G)$,
and $d=idim(R(G))=\sum\limits_{i=1}^k n_i$.
\end{theorem}
\proof By Algorithm \ref{A:FDL} and Theorem \ref{T:General}.
\qed\\

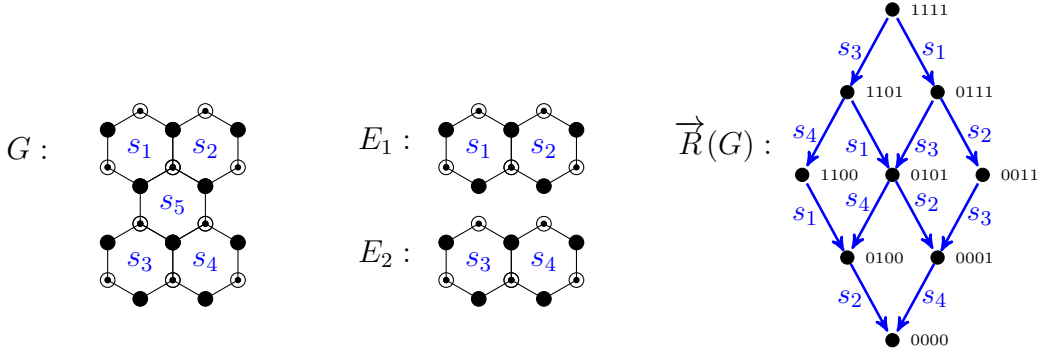
\begin{figure}[h!]
\begin{center}
\tikzset{Bullet/.style={draw=black, circle, minimum size=3pt, scale=0.5}}
\begin{tikzpicture}[scale=0.5]

  \begin{scope}[xshift=-3cm]
   \node [black]{$G:$};
  \end{scope}
  
  \draw [color=black, mark=*] plot[samples at={-150,-90,...,150,210},  variable=\x] 
  (\x:1);
    \node [Bullet] at (-150:1){};
  \node  [Bullet, fill] at (-90:1) {};
  \node [Bullet] at (-30:1) {} ;
  \node [Bullet, fill] at (30:1) {} ;
  \node  [Bullet] at (90:1) {};
  \node [Bullet, fill] at (150:1) {};
  \node [blue]{$s_1$};
  
   \begin{scope}[xshift=1.732 cm]
    \draw [color=black,mark=*] plot[samples at={-150,-90,...,210},variable=\x] 
  (\x:1);
      \node [Bullet] at (-150:1){};
  \node  [Bullet, fill] at (-90:1) {};
  \node [Bullet] at (-30:1) {} ;
  \node [Bullet, fill] at (30:1) {} ;
  \node  [Bullet] at (90:1) {};
  \node [Bullet, fill] at (150:1) {};
  \node[blue]{$s_2$};
  \end{scope}

 \begin{scope}[xshift=0.866 cm,  yshift=-1.5cm]
    \draw [color=black,mark=*] plot[samples at={-150,-90,...,210},variable=\x] 
  (\x:1);
      \node [Bullet] at (-150:1){};
  \node  [Bullet, fill] at (-90:1) {};
  \node [Bullet] at (-30:1) {} ;
  \node [Bullet, fill] at (30:1) {} ;
  \node  [Bullet] at (90:1) {};
  \node [Bullet, fill] at (150:1) {};
  \node[blue]{$s_5$};
  \end{scope}

  \begin{scope}[yshift=-3cm]
    \draw [color=black,mark=*] plot[samples at={-150,-90,...,210},variable=\x] 
  (\x:1);
      \node [Bullet] at (-150:1){};
  \node  [Bullet, fill] at (-90:1) {};
  \node [Bullet] at (-30:1) {} ;
  \node [Bullet, fill] at (30:1) {} ;
  \node  [Bullet] at (90:1) {};
  \node [Bullet, fill] at (150:1) {};
   \node [blue]{$s_3$};
  \end{scope}
  
    \begin{scope}[xshift=1.732 cm, yshift=-3cm]
    \draw [color=black,mark=*] plot[samples at={-150,-90,...,210},variable=\x] 
  (\x:1);
      \node [Bullet] at (-150:1){};
  \node  [Bullet, fill] at (-90:1) {};
  \node [Bullet] at (-30:1) {} ;
  \node [Bullet, fill] at (30:1) {} ;
  \node  [Bullet] at (90:1) {};
  \node [Bullet, fill] at (150:1) {};
   \node [blue]{$s_4$};
  \end{scope}

\begin{scope}[xshift=6.5cm, yshift=0.2cm]
\node[draw=none] {$E_1:$};
\end{scope}

   \begin{scope}[xshift=9 cm]
  \draw [color=black, mark=*] plot[samples at={-150,-90,...,150,210},  variable=\x] 
  (\x:1);
    \node [Bullet] at (-150:1){};
  \node  [Bullet, fill] at (-90:1) {};
  \node [Bullet] at (-30:1) {} ;
  \node [Bullet, fill] at (30:1) {} ;
  \node  [Bullet] at (90:1) {};
  \node [Bullet, fill] at (150:1) {};
  \node [blue]{$s_1$};
  \end{scope}
  
   \begin{scope}[xshift=10.732 cm]
    \draw [color=black,mark=*] plot[samples at={-150,-90,...,210},variable=\x] 
  (\x:1);
      \node [Bullet] at (-150:1){};
  \node  [Bullet, fill] at (-90:1) {};
  \node [Bullet] at (-30:1) {} ;
  \node [Bullet, fill] at (30:1) {} ;
  \node  [Bullet] at (90:1) {};
  \node [Bullet, fill] at (150:1) {};
  \node[blue]{$s_2$};
  \end{scope}

\begin{scope}[xshift=6.5cm, yshift=-2.8cm]
\node[draw=none] {$E_2:$};
\end{scope}
  \begin{scope}[xshift=9cm, yshift=-3cm]
    \draw [color=black,mark=*] plot[samples at={-150,-90,...,210},variable=\x] 
  (\x:1);
      \node [Bullet] at (-150:1){};
  \node  [Bullet, fill] at (-90:1) {};
  \node [Bullet] at (-30:1) {} ;
  \node [Bullet, fill] at (30:1) {} ;
  \node  [Bullet] at (90:1) {};
  \node [Bullet, fill] at (150:1) {};
   \node [blue]{$s_3$};
  \end{scope}
  
    \begin{scope}[xshift=10.732 cm, yshift=-3cm]
    \draw [color=black,mark=*] plot[samples at={-150,-90,...,210},variable=\x] 
  (\x:1);
      \node [Bullet] at (-150:1){};
  \node  [Bullet, fill] at (-90:1) {};
  \node [Bullet] at (-30:1) {} ;
  \node [Bullet, fill] at (30:1) {} ;
  \node  [Bullet] at (90:1) {};
  \node [Bullet, fill] at (150:1) {};
   \node [blue]{$s_4$};
  \end{scope}

\tikzset{>=stealth', shorten > = 1pt,  el/.style = {inner sep=2pt, align=left},every label/.append style = {font=\tiny}};

\begin{scope}[xshift=15.5cm, yshift=0.2cm]
\node[draw=none] {$\overrightarrow{R}(G):$};
\end{scope}

\begin{scope}[xshift=20cm, , yshift=-0.7cm]
\node [draw,shape=circle,  inner sep=0pt, minimum size=5pt, fill, label=right:$0101$] (1a) {};
\node [draw,shape=circle, inner sep=0pt, minimum size=5pt, fill, label=right:$1101$] (2a) [above of=1a,  xshift=-0.6cm, yshift=0.1cm] {};
\node [draw,shape=circle, inner sep=0pt, minimum size=5pt, fill, label=right:$0001$] (3a) [below of=1a, xshift=0.6cm, yshift=-0.1cm] {};
\node [draw,shape=circle, inner sep=0pt, minimum size=5pt, fill, label=right:$0100$] (4a) [below of=1a,  xshift=-0.6cm, yshift=-0.1cm] {};
\node [draw,shape=circle, inner sep=0pt, minimum size=5pt, fill, label=right:$0000$] (8a) [below of=1a, yshift=-1.2cm] {};
\node [draw,shape=circle, inner sep=0pt, minimum size=5pt, fill, label=right:$1100$] (5a) [left of=1a,  xshift=-0.2cm] {};
\node [draw,shape=circle, inner sep=0pt, minimum size=5pt, fill, label=right:$0111$] (6a) [above of=1a, xshift=0.6cm, yshift=0.1cm] {};
\node [draw,shape=circle, inner sep=0pt, minimum size=5pt, fill, label=right:$0011$] (9a) [right of=1a,  xshift=0.2cm] {};
\node [draw,shape=circle, inner sep=0pt, minimum size=5pt, fill, label=right:$1111$] (7a) [above of=1a, yshift=1.2cm] {};
\path[every node/.style={ font=\sffamily\small},   inner sep=0.5pt, line width=1.pt]
[<-] (1a) edge [color=blue] node[el, left, pos=0.3] { $s_1$} (2a)
[<-] (4a) edge [color=blue] node[el, left, pos=0.5] { $s_1$} (5a)
[<-] (6a) edge [color=blue] node[el, right, pos=0.5] { $s_1$} (7a)
[->] (4a) edge [color=blue] node[el, left, pos=0.5] { $s_2$} (8a)
[->] (1a) edge [color=blue] node[el, right, pos=0.3] { $s_2$} (3a)
[->] (6a) edge [color=blue] node[el, right, pos=0.5] { $s_2$} (9a)
[->] (1a) edge [color=blue] node[el, left, pos=0.3] { $s_4$} (4a)
[->] (2a) edge [color=blue] node[el, left, pos=0.5] { $s_4$} (5a)
[->] (3a) edge [color=blue] node[el, right, pos=0.5] { $s_4$} (8a)
[<-] (1a) edge [color=blue] node[el, right, pos=0.3] { $s_3$} (6a)
[<-] (2a) edge [color=blue] node[el, left, pos=0.5] { $s_3$} (7a)
[<-] (3a) edge [color=blue] node[el, right, pos=0.5] { $s_3$} (9a);
\end{scope}
\end{tikzpicture}
\end{center}
\caption{\label{Dimension-Example} A plane weakly elementary bipartite graph $G$ with two elementary components $E_1$ and $E_2$,
and a binary coding $[\phi]=(\phi^{[E_1]}, \phi^{[E_2]})$ of length $4$ on $\mathcal{M}(G)$.}
\end{figure}

\section{Remarks}
Our main results solve the fundamental cases of two open problems raised in \cite{ZLS08}.

\begin{remark}Theorem \ref{T:G-IsometricDimension} 
determines the smallest integer $m$ when the finite distributive lattice $\mathbb{M}(G)=(\mathcal{M}(G), \le_L)$ can be embedded into $\{0,1\}^m$.
\end{remark}

\begin{remark}
Theorem \ref{T:A-General}  provides an algorithm to generate a binary coding for all 1-factors of a plane bipartite graph $G$ which has the minimum 
$\mathrm{idim}(R(G))$ without generating $\mathcal{M}(G)$. Such plane bipartite graphs include all benzenoid graphs without coronenes 
as nice subgraphs.
\end{remark}

\begin{remark}
Theorem \ref{T:General} and Corollary \ref{C:Benzenoid}
provide the improved versions with $\mathrm{idim}(R(G))$ for Theorem 4.3 and Corollary 4.4 in \cite{ZLS08}, respectively.
\end{remark}

Theorem 4.3 in \cite{ZLS08} states that if $G$ is a plane weakly elementary bipartite graph with the set of finite faces $\mathcal{F}$,
then $\phi$ is an embedding of $\mathbb{M}(G)$ into $\{0,1\}^{\mathcal{F}}$ if and only if for any two disjoint cycles 
that form a nice subgraph their interiors are disjoint.
It is known that any benzenoid system is a plane weakly elementary bipartite graph \cite{ZZ00}, 
and an elementary benzenoid system  $G$ has no coronenes as its nice subgraphs 
if and only if for any pair of vertex disjoint cycles that form a nice subgraph of $G$ their interiors are disjoint \cite{ZC86}.
Corollary 4.4 in \cite{ZLS08} states that if $G$ is a benzenoid system, then
$\phi$ is an embedding of $\mathbb{M}(G)$ into $\{0,1\}^{\mathcal{F}}$ if and only if $G$ has no coronenes as its nice subgraphs.

We observe that $|\mathcal{F}|$ is not necessarily the isometric dimension of $R(G)$ in Theorem 4.3 and Corollary 4.4 in \cite{ZLS08}.
This can be seen by the example given in Figure \ref{Dimension-Example},
where the graph $G$ satisfies the properties stated in both Theorem 4.3 and Corollary 4.4.
So, $\phi$ is an embedding of $\mathbb{M}(G)$ into $\{0,1\}^{\mathcal{F}}$ where $|\mathcal{F}|=5$. 
But $\phi_M(s_5)=0$ for all $M \in \mathcal{M}(G)$, and by Theorem \ref{T:General}, the isometric dimension of $R(G)$ is $4$,
which is the summation of the numbers of finite faces of $E_1$ and $E_2$.
Figure \ref{Dimension-Example} provides an embedding $[\phi]$ of $\mathbb{M}(G)$ into $\{0,1\}^{4}$.

\vskip 0.2in
\noindent{\bf Acknowledgment:} 
  
 The research work is supported by the Research Development Grant (RDG) from Penn State University, Beaver Campus. 


\end{document}